\documentclass[a4paper,12pt]{article}

\usepackage[english]{babel} 
\usepackage{amsmath,amsfonts,amssymb,amsthm,bbm,graphics,epsfig,psfrag,epstopdf}
\usepackage[active]{srcltx}
\usepackage{paralist,subfigure,float}
\usepackage{color,soul}
\usepackage[dvipsnames]{xcolor}
\usepackage{hyperref} 
\setstcolor{red}

\newtheorem{theorem}{Theorem}
\newtheorem{corollary}[theorem]{Corollary}

\newtheorem{proposition}[theorem]{Proposition}

\newtheorem{remark}[theorem]{Remark}

\newtheorem{conjecture}[theorem]{Conjecture}
\newtheorem{hypothesis}[theorem]{Hypothesis}
\newtheorem{question}[theorem]{Question}

\def\N{\mathbb{N}}
\def\Z{\mathbb{Z}}

\def\R{\mathbb{R}}
\def\epsilon{\varepsilon}

\let\vp=\varphi

\let\t=\widetilde
\let\ol=\overline

\let\.=\cdot
\let\0=\emptyset
\let\mc=\mathcal

\def\O{\Omega}

\def\W{\mc{W}}

\def\tilde{\widetilde}

\DeclareMathOperator{\dist}{dist}

\def\thm#1{Theorem~\ref{thm:#1}}
\def\seq#1{(#1_n)_{n\in\N}}

\def\as{\quad\text{as }\;}
\def\inn{\quad\text{in }\;}

\def\tilde{\widetilde}
\def\1{\mathbbm{1}}
\def\Sph{\mathbb{S}^{N-1}}

\def\aos{asymptotic one-dimensional symmetry}

\newenvironment{formula}[1]{\begin{equation}\label{#1}}{\end{equation}\noindent}
\def\Fi#1{\begin{formula}{#1}}
\def\Ff{\end{formula}\noindent}

\newcommand{\be}{\begin{equation}}
\newcommand{\ee}{\end{equation}}
\newcommand{\baa}{\begin{array}}
\newcommand{\eaa}{\end{array}}
\newcommand{\ba}{\begin{eqnarray}}
\newcommand{\ea}{\end{eqnarray}}

\begin{document}
\date{}
\title{\bf{Spreading sets and one-dimensional symmetry for reaction-diffusion equations}}
\author{Fran\c cois Hamel$^{\hbox{\small{ a}}}$ and Luca Rossi$^{\hbox{\small{ b,c }}}$\thanks{This work has received funding from Excellence Initiative of Aix-Marseille Universit\'e~-~A*MIDEX, a French ``Investissements d'Avenir'' programme, and from the French ANR RESISTE (ANR-18-CE45-0019) project. The first author acknowledges support of the Institut Henri Poincaré (UAR 839 CNRS-Sorbonne Universit\'e), LabEx CARMIN (ANR-10-LABX-59-01), and Universit\`a degli Studi di Roma La Sapienza, where he was Sapienza Visiting Professor and where part of this work was done.}\\
\\
\footnotesize{$^{\hbox{a }}$Aix Marseille Univ, CNRS, I2M, Marseille, France}\\
\footnotesize{$^{\hbox{b }}$SAPIENZA Univ Roma, Istituto ``G.~Castelnuovo'', Roma, Italy}\\
\footnotesize{$^{\hbox{c }}$CNRS, EHESS, CAMS, Paris, France}}
\maketitle

\begin{center}
S\'eminaire Laurent Schwartz 2022
\end{center}

\begin{abstract}
\noindent{}We consider reaction-diffusion equations $\partial_tu=\Delta u+f(u)$ in the whole space $\R^N$ and we are interested in the large-time dynamics of solutions ranging in the interval $[0,1]$, with general unbounded initial support. 
Under the hypothesis of the existence of a traveling front connecting $0$ and $1$ with a positive speed, we discuss the existence of spreading speeds and spreading sets, which describe the large-time global shape of the level sets of the solutions. The spreading speed in any direction is expressed as a Freidlin-G\"artner type formula. This formula holds under general assumptions on the reaction and for solutions emanating from initial conditions with general unbounded support, whereas most of earlier results were concerned with more specific reactions and compactly supported or almost-planar initial conditions. We then investigate the local properties of the level sets at large time. Some flattening properties of the level sets of the solutions, if initially supported on subgraphs, will be presented. We also investigate the special case of asymptotically conical-shaped initial conditions. For Fisher-KPP equations, we state some asymptotic local one-dimensional and monotonicity symmetry properties for the elements of the $\Omega$-limit set of the solutions, in the spirit of a conjecture of De Giorgi for stationary solutions of Allen-Cahn equations. Lastly, we present some logarithmic-in-time estimates of the lag of the position of the solutions with respect to that of a planar front with minimal speed, for initial conditions which are supported on subgraphs with logarithmic growth at infinity. Some related conjectures and open problems are also listed.
\end{abstract}
	
		
\section{Framework and two main questions}\label{sec:2questions}

We consider solutions of the reaction-diffusion equation
\Fi{homo}
\partial_t u=\Delta u+f(u),\quad t>0,\ \ x\in\R^N,
\Ff
with $N\ge2$ and initial conditions $u_0$ having unbounded support. More precisely, the reaction term~$f:[0,1]\to\R$ is of class $C^1([0,1])$ with $f(0)=f(1)=0$, and the initial conditions $u_0$ are assumed to be characteristic functions $\1_U$ of sets~$U$, i.e.
\be\label{defu0}
u_0(x)=\left\{\baa{ll}1&\hbox{if }x\in U,\vspace{3pt}\\ 0&\hbox{if }x\in\R^N\!\setminus\!U,\eaa\right.
\ee
where the initial support $U$ is an unbounded measurable subset of $\R^N$ (we use the term ``initial support", with an abuse of notation, to refer to the set~$U$ in the definition of~$u_0$). The Cauchy problem is well posed and, given $u_0$, there is a unique bounded classical solution $u$ of~\eqref{homo} such that~$u(t,\cdot)\to u_0$ as~$t\to0^+$ in~$L^1_{loc}(\R^N)$. More general initial conditions $0\leq u_0\leq 1$ for which the upper level set $\{x\in\R^N:u_0(x)\geq\theta\}$ is at bounded Hausdorff distance from the support of $u_0$, where $\theta\in(0,1)$ is a suitable value depending on $f$, could be 
envisioned, at the expense of some further assumptions on the reaction term~$f$. For the sake of simplicity of the presentation, we focus on initial conditions $u_0$ of the type~\eqref{defu0}.

Due to diffusion, the solution $u$ of~\eqref{homo}-\eqref{defu0} is smooth at positive times and satisfies
$$0<u<1\ \hbox{ in }(0,+\infty)\times\R^N$$
from the strong parabolic maximum principle, provided the Lebesgue measures of $U$ and~$\R^N\setminus U$ are positive. However, from parabolic estimates, at each finite time, $u$ stays close to $1$ or $0$ in subregions of $U$ or $\R^N\setminus U$ which are far away from $\partial U$. One first goal is to describe the 
shape at large time of the regions where~$u$ stays close to $1$ or $0$. How do these regions move and possibly spread in any direction? A fundamental issue is to understand whether and how the solution keeps a memory at large time of its initial support~$U$.
 A~basic question is the following:

\begin{question}\label{q1}
For a given vector $e\in\R^N$ with Euclidean norm equal to $1$, is there a {\em spreading speed}~$w(e)$ such that
\be\label{ass0}\left\{\baa{lll}
u(t,cte)\to1 & \hbox{as $t\to+\infty$} & \hbox{for every $0\le c<w(e)$},\vspace{5pt}\\
u(t,cte)\to0 & \hbox{as $t\to+\infty$} & \hbox{for every $c>w(e)$}.\eaa\right.
\ee
Can one find a formula for $w(e)$ and, if any, how does $w(e)$ depend on~$e$ and the initial~support~$U$? Is there a uniformity with respect to~$e$ in~\eqref{ass0} and are there {\rm spreading sets} which describe the global shape of the level sets of~$u$ at large time?
\end{question}

Question~\ref{q1} will be addressed in Theorems~\ref{th1}-\ref{th4} below. Some conditions need to be imposed on the initial support~$U$, since otherwise the answer to Question~\ref{q1} can be negative in general (some counter-examples will be presented at the end of Section~\ref{sec:FG}). We also point out that, in~\eqref{ass0}, the speed $w(e)$ can possibly be $+\infty$ in some directions $e$, and this actually occurs in the directions around which $U$ is unbounded, 
in a sense that will be made precise~later.

Another question is concerned with the description of the profile of the solution around its level sets at large time. With this respect, we investigate two classes of properties: the {\em flattening} of the level sets, and the {\em asymptotic one-dimensional symmetry} of the solution. The latter is expressed in terms of the notion of $\Omega$-limit set, which is defined as follows: for a given bounded function $u:\R^+\times\R^N\to\R$, the set	
\be\label{def:Omega}\begin{array}{ll}
\Omega(u)\!:=\!\big\{&\!\!\!\!\!\psi\in L^\infty(\R^N)\ :\   u(t_n,x_n+\.)\to\psi\text{ in $L^\infty_{loc}(\R^N)$ as $n\to+\infty$,}\\
& \!\!\!\!\!\text{for some sequences $(t_n)_{n\in\N}$ in $\R^+$ diverging to $+\infty$ and $(x_n)_{n\in\N}$ in $\R^N$}\big\}
\end{array}
\ee
is called the $\O$-limit set of $u$. Roughly speaking, 
the $\O$-limit set contains all possible asymptotic profiles of the function as $t\to+\infty$. 
For any bounded solution $u$ of~\eqref{homo}, the set~$\Omega(u)$ is not empty and is included in $C^2(\R^N)$, 
by standard parabolic estimates. Motivated by some known results in the literature, the following question naturally arises.

\begin{question}\label{q2}
Let $u$ be a solution to~\eqref{homo} emerging from an initial datum $u_0=\1_U$. Is~it~true~that any function $\psi\in\Omega(u)$ is of the form
$$\psi(x)\equiv\Psi(x\.e),$$
for some~$e\in\Sph$ and $\Psi:\R\to\R$? If the answer to the question is positive, we then say that $u$ satisfies the {\rm asymptotic one-dimensional symmetry}.
\end{question}

For the answer to Question~\ref{q2} to possibly be affirmative, some conditions on $f$ and~$U$ need to be imposed, as shown by some counter-examples 
presented in Section~\ref{secaos}. We will also review in that section some known positive results which hold in the case where the initial support~$U$ is bounded, or when it is at finite Hausdorff distance from a half-space, under some assumptions on $f$. We will see how such results can be extended for a nonlinearity $f$ of the Fisher-KPP type, see condition~\eqref{fkpp} below, giving a positive answer to Question~\ref{q2} when~$U$ fulfills (in particular) a uniform interior ball condition and is convex, or, more generally, is at bounded Hausdorff distance from a convex set, see Theorem~\ref{thm:DG} below. These conditions on $U$ are actually a very particular instance of the geometric hypotheses under which we derive our most general result about the \aos, Theorem~\ref{thm:DGgeneral} below. Question~\ref{q2} reclaims the De Giorgi conjecture about solutions of the Allen-Cahn equation (that is, stationary solutions of the reaction-diffusion equation  $\Delta u+u(1-u)(u-1/2)=0$ in $\R^N$, obtained after a change of unknown from the original Allen-Cahn equation), see~\cite{dGconj}.

The situation considered here can be viewed as a  counterpart of many works devoted to the large-time dynamics of solutions of~\eqref{homo} with initial conditions~$u_0$ that are compactly supported or converge to $0$ at infinity. We refer to e.g.~\cite{AW,DM1,LK,MZ1,MZ2,Z} for extinction/invasion results in terms of the size and/or the amplitude of the initial condition~$u_0$ for various functions $f$, and to~\cite{DM1,DP,MP1,MP2,P1} for general local convergence and quasiconvergence results. For the invading solutions~$u$ (that is, those converging to $1$ locally uniformly in $\R^N$ as $t\to+\infty$) with localized initial conditions, further estimates on the location and shape at large time of the level sets have been established in~\cite{D,G,J,RRR,R2,R,U2}.

The case of general unbounded initial supports $U$ has been much less investigated in the lite\-rature. One immediately sees that, for general unbounded sets $U$, Question~\ref{q1} is much more intricate than in the case of bounded sets $U$, since the solutions $u$ can spread from all regions of the initial support $U$, that is, not only from a single bounded region. The sets~$U$ themselves can be bounded in some directions and unbounded in others. 


\section{Two main hypotheses}\label{sec:hyp}

In this section, we list some notations and hypotheses which are used in the various main results. 
The hypotheses are expressed in terms of the solutions of~\eqref{homo} with more general initial conditions than characteristic functions, or actually in terms of the reaction term $f$ solely. We then discuss the logical link between these hypotheses. We let ``$|\ |$'' and ``$\ \cdot\ $'' denote respectively the Euclidean norm and inner product in $\R^N$, $B_r(x):=\{y\in\R^N:|y-x|<r\}$ be the open Euclidean ball of center $x\in\R^N$ and radius $r>0$, $B_r:=B_r(0)$, and $\Sph:=\{e\in\R^N:|e|=1\}$ be the unit Euclidean sphere of~$\R^N$. The distance of a point $x\in\R^N$ from a set $A\subset\R^N$ is given by $\dist(x,A):=\inf\big\{|y-x|:y\in A\big\}$, with the convention $\dist(x,\emptyset)=+\infty$.  We also call $(\mathrm{e}_1,\cdots,\mathrm{e}_N)$ the canonical basis of $\R^N$, that is, $\mathrm{e}_i:=(0,\cdots,0,1,0,\cdots,0)$ for $1\le i\le N$, where $1$ is the $i$th coordinate of $\mathrm{e}_i$.

Since both $0$ and $1$ are steady states, the question of the interplay between these two states and the diffusion is intricate. One way to differentiate the roles of~$0$ and~$1$ is to assume that the state $1$ is more attractive than $0$, in the sense that it attracts the solutions of~\eqref{homo} --~not necessarily satisfying~\eqref{defu0}~-- that are ``large enough" in large balls at initial~time. 

\begin{hypothesis}\label{hyp:invasion}
The {\rm invasion property} occurs for any solution $u$ of~\eqref{homo} with a ``large enough'' initial datum $u_0$, that~is, there exist $\theta\in(0,1)$ and $\rho>0$ such that if
\be\label{hyptheta}
\theta\,\1_{B_\rho(x_0)}\le u_0\le1\inn\R^N,
\ee
for some $x_0\in\R^N$, then $u(t,x)\to1$ as $t\to+\infty$, locally uniformly with respect to $x\in\R^N$.
\end{hypothesis}
	
If $f$ satisfies the following conditions:
\Fi{HTconditions}
f>0 \text{ \ in $(0,1)$\quad and \quad }\liminf_{s\to0^+}\frac{f(s)}{s^{1+2/N}}>0,
\Ff 
then Hypothesis~\ref{hyp:invasion} is satisfied with any $\theta\in(0,1)$ and $\rho>0$, see~\cite{AW}; this property is known as the {\em hair trigger effect}. If $f>0$ in~$(0,1)$ (without any further assumption on the behavior of $f$ at $0^+$), then Hypothesis~\ref{hyp:invasion} is still satisfied with any $\theta\in(0,1)$, and with $\rho>0$ large enough. Hypothesis~\ref{hyp:invasion} holds as well if~$f$ is of the ignition type, that is,
\be\label{ignition}
\exists\,\alpha\in(0,1),\ \ f=0\hbox{ in $[0,\alpha]$ and $f>0$ in $(\alpha,1)$},
\ee
and $\theta$ in Hypothesis~\ref{hyp:invasion} can be any real number in the interval $(\alpha,1)$, provided $\rho>0$ is large enough. For a bistable function $f$ satisfying
\be\label{bistable}
\exists\,\alpha\in(0,1),\ \ f<0\hbox{ in }(0,\alpha)\hbox{ and }f>0\hbox{ in }(\alpha,1),
\ee
Hypothesis~\ref{hyp:invasion} is equivalent to $\int_0^1f(s)\,ds>0$, see~\cite{AW,FM}, and in that case $\theta$ in Hypo\-thesis~\ref{hyp:invasion} can be any real number in $(\alpha,1)$, provided $\rho>0$ is large enough. However, without the lower bound in condition~\eqref{hyptheta}, the solutions $u$ may not converge to $1$ at $t\to+\infty$ locally uniformly in $\R^N$, as easily seen for instance with functions $f$ of the types~\eqref{ignition} or~\eqref{bistable}, when $\|u_0\|_{L^1(\R^N)}$ is small enough. For a tristable function $f$ satisfying
\be\label{tristable}
\exists\,0<\alpha<\beta<\gamma<1,\quad f<0\hbox{ in $(0,\alpha)\cup(\beta,\gamma)\ $ and $\ f>0$ in $(\alpha,\beta)\cup(\gamma,1)$},
\ee
then it easily follows from~\cite{FM} that Hypothesis~\ref{hyp:invasion} is equivalent to the positivity of both integrals $\int_{\beta}^1f$ and $\int_0^1f$.

More generally speaking, it actually turns out that Hypothesis~\ref{hyp:invasion} is equivalent to the following two simple simultaneous conditions on the function~$f$, see~\cite{DP,P3}:
\be\label{hyp:theta}
\exists\,\theta\in(0,1),\quad f>0\hbox{ in }[\theta,1),
\ee
and
\be\label{intf}
\forall\,t\in[0,1),\quad\int_t^1f(s)\,ds>0.
\ee
Furthermore, $\theta$ can be chosen as the same real number in Hypothesis~$\ref{hyp:invasion}$ and in~\eqref{hyp:theta}. In particular, Hypothesis~\ref{hyp:invasion} is satisfied if $f\ge0$ in $[0,1]$ and if condition~\eqref{hyp:theta} holds. Notice however that condition~\eqref{hyp:theta} alone is not enough to guarantee Hypothesis~\ref{hyp:invasion}, since functions~$f$ of the type~\eqref{bistable} satisfy~\eqref{hyp:theta} but do not satisfy Hypothesis~\ref{hyp:invasion} if $\int_0^1f\le0$. Similarly, condition~\eqref{intf} alone is not enough to guarantee Hypothesis~\ref{hyp:invasion}, since there are $C^1([0,1])$ functions $f$ which vanish at $0$ and $1$ and satisfy~\eqref{intf} but not~\eqref{hyp:theta}: consider for instance $f$ defined by $f(1)=0$ and $f(s)=s(1-s)^3\sin^2(1/(1-s))$ for $s\in[0,1)$. Notice that, from the equivalence between Hypothesis~\ref{hyp:invasion} and~\eqref{hyp:theta}-\eqref{intf}, Hypothesis~\ref{hyp:invasion} is then independent of the dimension $N$, whereas, for a function $f$ which is positive in $(0,1)$, the hair trigger effect (that is, the arbitrariness of $\theta\in(0,1)$ and $\rho>0$ in Hypothesis~\ref{hyp:invasion}) depends on~$N$ (for instance, for the function $f(s)=s^p(1-s)$ with $p\ge1$, Hypothesis~\ref{hyp:invasion} holds in any dimension $N\ge1$, but the hair trigger effect holds if and only if $p\le1+2/N$, see~\cite{AW}).

In the large time dynamics of the solutions of the Cauchy problem~\eqref{homo}, a crucial role is played by the {\em traveling front} solutions connecting the steady states $1$ and $0$, defined as
$$u(t,x)=\vp(x\cdot e-ct)$$
with $c\in\R$, $e\in\Sph$, and
\be\label{limitsvp}
1=\vp(-\infty)>\vp(z)>\vp(+\infty)=0\quad\hbox{for all $z\in\R$}.
\ee
The level sets of these solutions are hyperplanes orthogonal to $e$, traveling with constant speed $c$ in the direction~$e$. If any, the profile $\vp$ is necessarily decreasing and unique up to shifts, for a given~$c$. Most of the main results are derived under the following hypothesis:

\begin{hypothesis}\label{hyp:minimalspeed}
Equation~\eqref{homo} admits a traveling front connecting $1$ to~$0$ with speed $c_0>0$.
\end{hypothesis}
	
Hypothesis~\ref{hyp:minimalspeed} is fulfilled for instance if $f>0$ in $(0,1)$, or if $f$ is of the ignition type~\eqref{ignition}, or if $f$ is of the bistable type~\eqref{bistable} with $\int_0^1f(s)\,ds>0$ (in the last two cases, the speed $c_0$ is unique), see~\cite{AW,FM,KPP}. Hypothesis~\ref{hyp:minimalspeed} is also satisfied for some functions~$f$ having multiple oscillations in the interval $[0,1]$, see Remark~\ref{rem6} below. It actually turns out that Hypothesis~\ref{hyp:minimalspeed} is equivalent to the existence of a positive minimal speed $c^*$ of traveling fronts connecting $1$ to $0$, and that Hypothesis~\ref{hyp:minimalspeed} also implies Hypo\-thesis~\ref{hyp:invasion} and further spreading properties for the solutions of~\eqref{homo} fulfilling the conditions of Hypothesis~\ref{hyp:invasion}:

\begin{proposition}\label{pro:spreadingFL}{\rm{\cite{HR1}}}
Assume Hypothesis~$\ref{hyp:minimalspeed}$. Then equation~\eqref{homo} admits a traveling front connecting $1$ to $0$ with minimal speed $c^*$, and $c^*>0$. Furthermore, Hypothesis~$\ref{hyp:invasion}$ is fulfilled and, for any solution $u$ as in Hypothesis~$\ref{hyp:invasion}$, 
it holds that
\be\label{c<c*}
\forall\,c\in[0,c^*),\quad\min_{|x|\le ct}u(t,x)\to1\as t\to+\infty.
\ee
Lastly, for any compactly supported initial datum $0\le u_0\le 1$, the solution $u$ of~\eqref{homo} satisfies
$$\forall\,c>c^*,\quad\sup_{|x|\geq ct}u(t,x)\to0\as t\to+\infty.$$
\end{proposition}

The minimality of $c^*$ means that~\eqref{homo} in dimension $N=1$ admits a solution of the form $\vp(x-c^*t)$ satisfying~\eqref{limitsvp}, 
and it does not admit a solution of the same type with~$c^*$ replaced by a smaller quantity (notice that, necessarily, $c^*\le c_0$ under the notation of Hypothesis~\ref{hyp:minimalspeed}). Proposition~\ref{pro:spreadingFL} answers affirmatively to Question~\ref{q1} under Hypothesis~$\ref{hyp:minimalspeed}$, in the 
very special case of compactly supported initial data satisfying~\eqref{hyptheta}, with $w(e)=c^*$ for all $e\in\Sph$. This can be viewed as a natural extension of some results of the seminal paper~\cite{AW}, which were originally obtained under more specific assumptions on $f$.

\begin{remark}\label{rem6}{\rm
Whereas Proposition~\ref{pro:spreadingFL} shows the implication ``Hypothesis~\ref{hyp:minimalspeed} $\Longrightarrow$ Hypothesis~\ref{hyp:invasion}", the converse implication is false in general. For instance, consider equation~\eqref{homo} with~$f$ satisfying~\eqref{tristable} together with $\int_0^\beta f>0$ and $\int_\beta^1f>0$, and let $c_1$ and $c_2$ be the unique (positive) speeds of the traveling fronts $\vp_1(x-c_1t)$ and $\vp_2(x-c_2t)$ connecting $\beta$~to~$0$ on the one hand, and~$1$ to~$\beta$ on the other hand. It follows from~\cite{FM} that, if~$c_1\ge c_2$, then Hypothesis~\ref{hyp:minimalspeed} is not satisfied, whereas Hypothesis~\ref{hyp:invasion} is. In that case, it turns out that the compactly supported initial conditions $u_0$ giving rise to invading solutions $u$ develop into a terrace of two expanding fronts with speeds~$c_1$ and~$c_2$, in the sense that $\inf_{B_{ct}}u(t,\cdot)\to1$ as $t\to+\infty$ if $0<c<c_2$ (resp. $\sup_{B_{c''t}\setminus B_{c't}}|u(t,\cdot)-\beta|\to0$ as $t\to+\infty$ if $c_2<c'<c''<c_1$, resp. $\sup_{\R^N\setminus B_{ct}}u(t,\cdot)\to0$ as $t\to+\infty$ if $c>c_1$). We refer to~\cite{DM2,DGM,GR,P2,P3} for further results on terraces in more general frameworks. On the other hand, still with~\eqref{tristable} and the positivity of  $\int_0^\beta f$ and $\int_\beta^1f$, Hypothesis~\ref{hyp:minimalspeed} is satisfied if (and, then, only~if) $c_1<c_2$, see~\cite{FM}.}
\end{remark}

	
\section{Spreading speeds and spreading sets}\label{sec:FG}

In this section, under Hypothesis~\ref{hyp:minimalspeed}, we investigate the notions of asymptotic spreading speeds and spreading sets for the solutions $u$ of~\eqref{homo}-\eqref{defu0} with general unbounded sets $U$ containing large enough balls. Such solutions $u$ then converge to $1$ as $t\to+\infty$ locally uniformly in $\R^N$, and even satisfy~\eqref{c<c*}, with $c^*>0$ given by Proposition~\ref{pro:spreadingFL}. We now want to provide a more precise description of the invasion of the state~$0$ by the state $1$. The invasion cannot be uniform in all directions in general, since it shall strongly depend on the initial support $U$. For~$e\in\Sph$, we then look for a quantity $w(e)\in(0,+\infty]$ satisfying~\eqref{ass0}, referred to as the {\em spreading speed} and representing the asymptotic speed at which the level sets between $0$ and $1$ move 
along the direction $e$. If any, it satisfies $w(e)\ge c^*$ by Proposition~\ref{pro:spreadingFL}. However, contrary to the case of compactly supported initial data satisfying~\eqref{hyptheta}, the spreading speed may not exist when the support of the initial condition is unbounded, see the comments at the end of this section.

Let us first introduce the notions of sets of directions ``around which~$U$ is bounded'' and ``around which~$U$ is unbounded'', for short the sets of {\it bounded directions} and of {\it unbounded directions}, defined by:
$$\mc{B}(U):=\Big\{\xi\in\Sph:\liminf_{\tau\to+\infty}\frac{\dist(\tau\xi,U)}{\tau}>0\Big\}$$
and
$$\mc{U}(U):=\Big\{\xi\in\Sph:\lim_{\tau\to+\infty}\frac{\dist(\tau\xi,U)}{\tau}=0\Big\}.$$
The sets~$\mc{B}(U)$ and $\mc{U}(U)$ are respectively open and closed relatively to $\Sph$. 
The condition $\xi\in\mc{B}(U)$ is equivalent to the existence of an open cone $\mc{C}$ containing the 
ray~$\R^+\{\xi\}=\{\tau\,\xi:\tau>0\}$ such that $U\cap\mc{C}$ is bounded. Conversely, for any $\xi\in\mc{U}(U)$ and any open cone $\mc{C}$ containing the ray~$\R^+\{\xi\}$, the set $U\cap\mc{C}$ is unbounded. We also define the notion of positive-distance-interior~$U_\rho$ (with $\rho>0$) of the set $U$ as
$$U_\rho:=\big\{x\in U:\dist(x,\partial U)\ge\rho\big\}.$$

The first main result shows the existence of and a formula for the spreading speeds, providing a positive answer to the first part of Question~\ref{q1}.

\begin{theorem}\label{th1}{\rm{\cite{HR1}}}
Assume that Hypothesis $\ref{hyp:minimalspeed}$ holds, let~$c^*>0$ and $\rho>0$ be given by Proposition~$\ref{pro:spreadingFL}$ and Hypothesis~$\ref{hyp:invasion}$, and let~$u$ be the solution of~\eqref{homo}-\eqref{defu0}, with $U_\rho\neq\emptyset$~and
\Fi{hyp:U}
\mc{B}(U)\cup\,\mc{U}(U_\rho)=\Sph.
\Ff
Then, for every $e\in\Sph$, there exists $w(e)\in[c^*,+\infty]$ such that~\eqref{ass0} holds, and even
\Fi{ass}\left\{\baa{lll}
\displaystyle\lim_{t\to+\infty}\,\Big(\min_{0\le s\le c}u(t, ste)\Big)\!&\! =1 & \text{for every }0\le c<w(e),\vspace{3pt}\\
\displaystyle\lim_{t\to+\infty}\,\Big(\sup_{s\ge c}u(t,ste)\Big)\!&\!=0 & \text{for every }c>w(e).\eaa\right.
\Ff
Furthermore, $w(e)$ is given explicitly by the variational formula
\Fi{FGgeneral}
w(e)=\sup_{\xi\in\mc{U}(U),\ \xi\.e\ge0}\ \frac{c^*}{\sqrt{1-(\xi\.e)^2}},
\Ff
with $w(e)=c^*$ if there is no $\xi\in\mc{U}(U)$ such that $\xi\.e\ge0$, and $w(e)=+\infty$ if $e\in\mc{U}(U)$.
\end{theorem}

Since $\mc{U}(U)$ is closed in $\Sph$, it follows from~\eqref{FGgeneral} and the above conventions that
$$\left\{\baa{ll}
w(e)=+\infty & \!\!\hbox{if and only if }e\in\mc{U}(U),\vspace{3pt}\\
w(e)>c^* & \!\!\hbox{if and only if there is $\xi\in\mc{U}(U)$ such that $\xi\cdot e>0$},\vspace{3pt}\\
w(e)=c^* & \!\!\hbox{if and only if there is no $\xi\in\mc{U}(U)$ such that $\xi\cdot e>0$}.\eaa\right.$$
For a set $U$ satisfying $\mc{U}(U)\neq\emptyset$, formula~\eqref{FGgeneral} can be rephrased in a more geometric way:
\Fi{FGgeometric}
w(e)\,=\,\frac{c^*}{\dist(e,\R^+\, \mc{U}(U))}\,=\,\frac{c^*}{\sin\vartheta}\,,
\Ff
where $\vartheta\in[0,\pi/2]$ is the minimum between~$\pi/2$ and the smallest angle between~$e$ and the directions in~$\mc{U}(U)$ (with the convention $c^*/0=+\infty$). 
This formula immediately reveals that the map $e\mapsto w(e)\in[c^*,+\infty]$ is continuous in $\Sph$. If $\mc{U}(U)=\emptyset$, then $w(e)=c^*$ for all $e\in\Sph$. If $U$ is bounded, condition~\eqref{hyp:U} is automatically satisfied and, if $U_\rho\neq\emptyset$, then~\eqref{ass} holds with $w(e)=c^*$ for all $e\in\Sph$, in agreement with Proposition~\ref{pro:spreadingFL}.

Formula~\eqref{FGgeneral} is called a Freidlin-G\"artner type formula, since Freidlin and G\"artner~\cite{FG} were the first to characterize the spreading speeds of solutions of reaction-diffusion equations in~$\R^N$ by a variational formula. 
They were actually concerned with spreading speeds for solutions of $x$-dependent reaction-diffusion equations 
of the Fisher-KPP type \cite{F,KPP} (for which $0<f(x,u)/u\le\frac{\partial f}{\partial u}(x,0)$ for all $(x,u)\in\R^N\times(0,1)$) 
with $f(x,u)$ periodic with respect to~$x$.
More precisely, it follows from~\cite{FG}, together with~\cite{BHN,BHN1,W}, that~\eqref{ass0} holds for these solutions, with
\be\label{formuleFG}
w(e)=\inf_{\xi\in\Sph,\,\xi\cdot e>0}\ \frac{c^*(\xi)}{\xi\cdot e}
\ee
for any $e\in\Sph$, where $c^*(\xi)$ denotes the minimal speed of pulsating fronts connecting $1$ to $0$ in the direction $\xi$ (a pulsating front connecting $1$ to $0$ with speed $c$ in the direction $\xi$ is a solution $u:\R\times\R^N\to(0,1)$ such that $u(t,x)=\phi(x\cdot\xi-ct,x)$, where $\phi(-\infty,x)=1$, $\phi(+\infty,x)=0$ uniformly in $x\in\R^N$, and $\phi$ has the same periodicity with respect to its second argument as the function $f$ or other coefficients of the equation, see e.g.~\cite{BH,SKT,W,X1,X2}). Such formulas for the spreading speeds of solutions with compactly supported initial conditions 
have been recently extended to more general reaction terms in~\cite{R1}. For reaction-diffusion equations with spatially periodic coefficients, the spreading speed $w(e)$ may depend on the direction~$e$, even for solutions with compactly supported initial conditions~$u_0$. However, the continuity of the map~$e\mapsto w(e)$ still holds for monostable, ignition or bistable reactions~$f$, as follows from~\cite{FG,R1} and from the (semi-)continuity of the minimal or unique speeds of pulsating traveling fronts with respect to the direction, see~\cite{AG,G,G2,R1} (but the continuity of the spreading speeds and even their existence do not hold in general when pulsating fronts connecting $1$ to $0$ do not exist anymore, see~\cite{GR}).

Remember that Hypothesis~\ref{hyp:minimalspeed} holds if $f>0$ in $(0,1)$, in the ignition case~\eqref{ignition}, and in the bistable case~\eqref{bistable} with~$\int_0^1f(s)ds>0$. In these cases, Theorem~\ref{th1} yields the existence of the spreading speeds satisfying~\eqref{ass0} and~\eqref{ass}, given by~\eqref{FGgeneral} as soon as the initial datum~$u_0=\1_U$ is associated with a set~$U\subset\R^N$ satisfying $U_\rho\neq\emptyset$ and~\eqref{hyp:U}. Moreover, in the case of a positive nonlinearity satisfying~\eqref{HTconditions}, $\rho>0$ can be arbitrarily small. On the other hand, the conclusions of Theorem~\ref{th1} do not hold in general without Hypothesis~\ref{hyp:minimalspeed}: as in Remark~\ref{rem6}, for a function $f$ of the type~\eqref{tristable} with $c_1>c_2$ (where~$c_1$ and~$c_2$ are the positive speeds of the traveling fronts $\vp_1(x-c_1t)$ and $\vp_2(x-c_2t)$ connecting $\beta$~to~$0$, and~$1$ to~$\beta$, respectively), the solutions $u$ of~\eqref{homo}-\eqref{defu0} with $U$ bounded (hence,~\eqref{hyp:U} is satisfied) and $U_\rho\neq\emptyset$ develop into a terrace of expanding fronts, ruling out the existence of $w(e)$ satisfying~\eqref{ass0}.

The geometric assumption~\eqref{hyp:U} is invariant under rigid transformations of~$U$. It holds for instance if $U$ is star-shaped and stays at a finite distance from its $\rho$-interior $U_\rho$, by~\cite{HR1}. 

The next result states the uniformity of~\eqref{ass0} with respect to the directions $e\in\Sph$, making more precise the answer to Question~\ref{q1}.

\begin{theorem}\label{th2}{\rm{\cite{HR1}}}
Under the assumptions of Theorem~$\ref{th1}$, for any compact set $C\subset\R^N$,
\Fi{ass-cpt}\left\{\baa{lll}
\displaystyle\lim_{t\to+\infty}\,\Big(\min_{x\in C}u(t, tx)\Big)\!&\! =1 & \text{if }\;C\subset\mc{W},\vspace{3pt}\\
\displaystyle\lim_{t\to+\infty}\,\Big(\max_{x\in C}u(t,tx)\Big)\!&\!=0 & \text{if }\;C\subset\R^N\setminus\ol{\mc{W}},\eaa\right.
\Ff
where $\mc{W}$ is the envelop set of the speeds $w(e)$'s, that is,
\be\label{asspre}
\W:=\big\{re\,:\,e\in \Sph,\ \ 0\leq r<w(e)\big\}.
\ee
\end{theorem}

Formula~\eqref{FGgeometric} reveals that $\mc{W}$ has the following simple geometric expression:
\be\label{asspre-formula}
\W=\R^+\, \mc{U}(U)\,+\,B_{c^*}
\ee
(with the convention that $\R^+\emptyset\,+\,B_{c^*}=B_{c^*}$ if $\mc{U}(U)=\emptyset$). 
Indeed, on the one hand, if $\mc{U}(U)=\emptyset$, then $w(e)\equiv c^*$ and~$\mc{W}=B_{c^*}$. 
On the other hand, if $\mc{U}(U)\neq\emptyset$, for any $e$ and $r\ge0$, one has
 $\dist(re,\R^+\mc{U}(U))=r\,\dist(e,\R^+\mc{U}(U))=rc^*/w(e)$ by~\eqref{FGgeometric} 
 (using the convention $rc^*/(+\infty)=0$), and therefore the equivalence between~\eqref{asspre} and~\eqref{asspre-formula} 
 follows. Formula~\eqref{asspre-formula} means that $\mc{W}$ is the $c^*$-neighborhood of 
 the positive cone generated by the directions $\mc{U}(U)$. 
 It implies that $\mc{W}$ is an open set which is either unbounded (when $\mc{U}(U)\neq\emptyset$), 
 or it coincides with $B_{c^*}$. For periodic Fisher-KPP equations, 
 formula~\eqref{formuleFG} for the spreading speeds of solutions with compactly supported initial conditions 
 means that the closure of the set $\W$ defined by~\eqref{asspre} coincides with the Wulff shape of the envelop set of the minimal speeds $c^*(\xi)$ of pulsating fronts and, since the map $\xi\mapsto c^*(\xi)\in(0,+\infty)$ is continuous in $\Sph$, the set $\overline{\W}$ would therefore be a convex compact set. For our pro\-blem~\eqref{homo}-\eqref{defu0}, the set $\W$ defined in~\eqref{asspre}-\eqref{asspre-formula} is not bounded as soon as $\mc{U}(U)\neq\emptyset$. Furthermore, it is not convex in general. For instance, if~$U\neq\emptyset$ is a non-convex closed cone, say with vertex $0$, then~$\R^+\mc{U}(U)\cup\{0\}=U$ and, from~\eqref{asspre-formula},~$\W$ is not convex either. Nevertheless, if $U$ is a general convex set, then $\R^+\mc{U}(U)\cup\{0\}$ and $\W$ are convex, from~\eqref{asspre-formula} again. More generally speaking, if there is a convex set $U'$ which lies at a finite Hausdorff distance from~$U$, then $\mc{U}(U)=\mc{U}(U')$ and therefore~$\W$ is convex, even if~$U$ itself is not.

Having in mind~\eqref{ass-cpt}, $\W$ is called a {\em spreading set} for~\eqref{homo}-\eqref{defu0}. We point out that~\eqref{ass-cpt} is stronger than~\eqref{ass0}, owing to the continuity of the map $e\mapsto w(e)$ in~$\Sph$. It also yields the first line of~\eqref{ass}. Compared to the first lines of~\eqref{ass0} and~\eqref{ass}, the first line of~\eqref{ass-cpt} provides an additional uniformity with respect to the directions~$e$. It also follows from~\eqref{ass-cpt} and the continuity of the map $e\mapsto w(e)\in[c^*,+\infty]$ that, for any $\sigma\in(0,1)$ and $A>0$,
$$\min_{x\,\in\,\sigma\overline{\W}\cap\overline{B_A}}u(t,tx)\to1\ \hbox{ as }t\to+\infty.$$
Formulas similar to~\eqref{ass-cpt} have been established for the solutions of more general heterogeneous equations or systems with compactly supported initial conditions and Fisher-KPP reactions~\cite{BES,BN,ES,LM,W}, 
bistable reactions~\cite{X2}, or even more general terms~\cite{R1,W}. The main difference is that, in these references, the spreading speeds and sets are bounded, unlike the spreading set $\W$ defined in~\eqref{asspre}-\eqref{asspre-formula}, which is unbounded as soon as $\mc{U}(U)\neq\emptyset$.

Theorems~\ref{th1}-\ref{th2} are concerned with the convergence towards $1$ and $0$ as $t\to+\infty$ along some rays or some dilated sets. The next two results provide a description of the asymptotic shape of the upper level sets of a solution $u$, defined for $\lambda\in(0,1)$ and $t>0$ by
\be\label{defElambda}
E_\lambda(t):=\big\{x\in\R^N:u(t,x)>\lambda\big\}.
\ee
That description involves the Hausdorff distance between some sets depending on~$E_\lambda(t)$ and~$t\W$. The Hausdorff distance is defined, for any pair of subsets $A,B\subset\R^N$, by  
$$d_{\mc{H}}(A,B):=\max\Big(\sup_{x\in A}\dist(x,B),\,\sup_{y\in B}\dist(y,A)\Big),$$
with the conventions that $d_{\mc{H}}(A,\emptyset)=d_{\mc{H}}(\emptyset,A)=+\infty$ if $A\neq\emptyset$ and $d_{\mc{H}}(\emptyset,\emptyset)=0$.

\begin{theorem}\label{th3}{\rm{\cite{HR1}}}
Under the assumptions of Theorems~$\ref{th1}$-$\ref{th2}$, it holds that
\Fi{Hloc}
\forall\,R>0,\ \ \forall\,\lambda\in(0,1),\ \ d_{\mc{H}}\Big(\overline{B_R}\cap\frac1t\, E_\lambda(t)\,,\,\overline{B_R}\cap \mc{W}\Big)\mathop{\longrightarrow}_{t\to+\infty}0.
\Ff
\end{theorem}

Theorem~\ref{th3} gives the approximation of $t^{-1}E_\lambda(t)$ by $\mc{W}$ as $t\to+\infty$, {\it locally} with respect to the Hausdorff distance. But the convergence is not {\it global} in general, and $d_{\mc{H}}(t^{-1}E_\lambda(t),\mc{W})\not\to0$ as $t\to+\infty$ in general, see the comments at the end of this section. However, it is global if $U$ is bounded and $U_{\rho}\neq\emptyset$, by Proposition~\ref{pro:spreadingFL}, with $\mc{W}=B_{c^*}$ in this case.

The following result provides an asymptotic global approximation of~$t^{-1}E_\lambda(t)$ by the family of sets $t^{-1}U+B_{c^*}$, under 
a different assumption on $U$.

\begin{theorem}\label{th4}{\rm{\cite{HR1}}}
Assume that Hypothesis~$\ref{hyp:minimalspeed}$ holds, let~$c^*>0$ and $\rho>0$ be given by Proposition~$\ref{pro:spreadingFL}$ and Hypothesis~$\ref{hyp:invasion}$, and let~$u$ be the solution of~\eqref{homo}-\eqref{defu0}, with $U_\rho\neq\emptyset$~and
\Fi{dUrho}
d_{\mc{H}}(U,U_\rho)<+\infty.
\Ff
Then, $d_{\mc{H}}(E_\lambda(t),U+B_{c^*t})=o(t)$ as $t\to+\infty$ for every $\lambda\in(0,1)$, that is,
\Fi{dH}
\forall\,\lambda\in(0,1),\ \ d_{\mc{H}}\Big(\frac{1}{t}\,E_\lambda(t) \,,\, \frac{1}{t}\,U+B_{c^*}\big)\mathop{\longrightarrow}_{t\to+\infty}0.
\Ff
\end{theorem}

Property~\eqref{dH} means that $E_\lambda(t)$ behaves at large time $t$ as the set $U$ thickened by~$c^*t$. A~sufficient condition for~\eqref{dUrho} to hold is that the set $U$ fulfills the uniform interior sphere condition of radius~$\rho$: in such case $d_{\mc{H}}(U,U_\rho)\leq2\rho$. In particular, if~$f$ satisfies condition~\eqref{HTconditions}, then Theorem~\ref{th4} applies to any  non-empty set $U$ which is uniformly $C^{1,1}$.  

We point out that a single formula like~\eqref{dH} valid for all $\lambda\in(0,1)$ does not hold in general without Hypo\-thesis~\ref{hyp:minimalspeed}. For instance, as in Remark~\ref{rem6}, consider a tristable function $f$ of the type~\eqref{tristable} with $c_1>c_2$ (where~$c_1$ and~$c_2$ are the positive speeds of the traveling fronts $\vp_1(x-c_1t)$ and $\vp_2(x-c_2t)$ connecting~$\beta$~to~$0$, and~$1$ to~$\beta$, respectively). Then, as follows from~\cite{DM2,DGM,FM,P2}, the solutions~$u$ of~\eqref{homo}-\eqref{defu0} with $U$ bounded and $U_\rho\neq\emptyset$ (hence,~\eqref{dUrho} is satisfied) are such that $d_{\mc{H}}(E_\lambda(t),U+B_{c_2t})=o(t)$ as $t\to+\infty$ if $\beta<\lambda<1$, respectively $d_{\mc{H}}(E_\lambda(t),U+B_{c_1t})=o(t)$ as $t\to+\infty$ if $0<\lambda<\beta$.

In order to enlighten our results stated above, we consider the important class of unbounded sets $U$
given by subgraphs of some functions (see also the next section for further results). Namely, we consider
\be\label{Ugamma}
U:=\big\{x=(x',x_N)\in\R^{N-1}\times\R:x_N\le\gamma(x')\big\},
\ee
with $\gamma:\R^{N-1}\to\R$ belonging to $L^\infty_{loc}(\R^{N-1})$. Assume for instance that $\gamma$ is of the form
$$\gamma(x')=\alpha\,|x'|+o(|x'|) \quad\text{as }\;|x'|\to+\infty,$$
for some $\alpha\in\R$. We see that $U_\rho\neq\emptyset$ for any $\rho>0$ and that 
$\mc{B}(U)=\big\{e\in\Sph:e_N>\alpha|e'|\big\}$ and $\mc{U}(U)=\mc{U}(U_\rho)=\big\{e\in\Sph:e_N\leq\alpha|e'|\big\}$. 
Thus~\eqref{hyp:U} is fulfilled, hence~\eqref{ass0},~\eqref{ass},~\eqref{ass-cpt} and~\eqref{Hloc} hold under 
Hypothesis~\ref{hyp:minimalspeed} on~$f$, by Theorems~\ref{th1}-\ref{th3}. 
However, the shape of the envelop set $\mc{W}$ given by~\eqref{asspre}-\eqref{asspre-formula} strongly depends on the sign of $\alpha$. If $\alpha>0$, then $\mc{W}=\{x\in\R^N:x_N< \alpha\,|x'|+c^*\sqrt{1+\alpha^2}\}$: it is a translation of the interior of the cone~$\R^+ \mc{U}(U)$, hence it is non-convex and not $C^1$. If~$\alpha<0$ then~$\mc{W}$ is still given by the $c^*$-neighborhood of the same cone $\R^+\mc{U}(U)$, which now becomes ``rounded'' in its upper part; indeed in such a case~$w(e)=c^*$ if~$e_N\ge|e'|/|\alpha|$, and $\mc{W}$ is convex and~$C^1$ (but not~$C^2$). If~$\alpha=0$ (which includes the case $\gamma$ bounded) then $\mc{W}=\{x\in\R^N:x_N<c^*\}$ is a half-space, with $w(e)=+\infty$ if~$e_N\le0$, and~$w(e)=c^*/e_N$ if~$e_N>0$. 
If $\gamma$ in~\eqref{Ugamma} satisfies
$$\frac{\gamma(x')}{|x'|}\to-\infty\hbox{ as $|x'|\to+\infty$},$$
then $\mc{B}(U)\!=\!\Sph\setminus\{-\mathrm{e}_N\}$ and $\mc{U}(U)\!=\!\mc{U}(U_\rho)\!=\!\{-\mathrm{e}_N\}$, with $\mathrm{e}_N=(0,\cdots,0,1)$. Hence~\eqref{hyp:U} is fulfilled and therefore, under Hypothesis~\ref{hyp:minimalspeed}, properties~\eqref{ass0},~\eqref{ass},~\eqref{ass-cpt} and~\eqref{Hloc} hold by Theorems~\ref{th1}-\ref{th3}, with $\mc{W}=-\R^+\mathrm{e}_N+B_{c^*}=\big\{x\in\R^N:|x'|<c^*,\ x_N\leq0\big\}\cup B_{c^*}$. This is a cylinder with a ``rounded'' top, which is convex and $C^1$, but not $C^2$. Lastly, if~$\gamma$ has uniformly bounded local oscillations, that is, if there is $M>0$ such that $|\gamma(x')-\gamma(y')|\le M$ for all $x'\in\R^{N-1},\,y'\in\R^{N-1}$ with $|x'-y'|\le1$ (this is the case if $\gamma$ is globally Lipschitz-continous), then condition~\eqref{dUrho} is fulfilled and~\eqref{dH} holds, by Theorem~\ref{th4}.

To complete this section, we present a list of situations where one or both hypotheses~\eqref{hyp:U} and~\eqref{dUrho} of Theorems~\ref{th1}-\ref{th4} do not hold and the conclusions~\eqref{ass0},~\eqref{ass},~\eqref{ass-cpt}, \eqref{Hloc} and~\eqref{dH} fail. The examples also show that the conditions~\eqref{hyp:U} and~\eqref{dUrho} on~$U$ cannot be compared in general. We also discuss the validity of the following convergences:
\Fi{Hausdorff}
\lim_{t\to+\infty}\frac1t\, E_\lambda(t)=\mc{W}=\lim_{t\to+\infty}\frac1t\,U+B_{c^*},
\Ff
that one may expect to hold but that actually fail in general.  The convergences~\eqref{Hausdorff} would be understood with respect to the Hausdorff distance (which does not guarantee the uniqueness of the limit). Notice first that, if~\eqref{hyp:U} is fulfilled together with $U_\rho\neq\emptyset$ and Hypothesis~\ref{hyp:minimalspeed}, then~\eqref{Hloc} holds and the limit of $t^{-1}E_\lambda(t)$, if any, must be the set $\mc{W}$ (that is, the Hausdorff distance between the limit set and $\mc{W}$ must be~$0$). All of the following instances refer to the equation~\eqref{homo} with logistic term $f(u)=u(1-u)$, for which Hypothesis~\ref{hyp:minimalspeed} holds, as well as the hair trigger effect, i.e., $\theta\in(0,1)$ and $\rho>0$ can be arbitrary in Hypothesis~\ref{hyp:invasion} and Proposition~\ref{pro:spreadingFL}, in which $c^*=2$. The sets $U$ listed below have non-empty interiors, and~\eqref{hyp:U} and~\eqref{dUrho} are understood here with~$\rho>0$ arbitrarily~small.

\begin{itemize}
\item The set $U:=\bigcup_{n\in\N}\overline{B_{2^n+1}}\setminus B_{2^n-1}$ fulfills~\eqref{dUrho} (hence~\eqref{dH} holds), but it violates~\eqref{hyp:U}, and~\eqref{ass0}, \eqref{ass},~\eqref{ass-cpt} and~\eqref{Hloc} all fail, for any function $w:\Sph\to[0,+\infty]$ and any star-shaped, open set~$\mc{W}\subset\R^N$, and moreover both limits in~\eqref{Hausdorff} do not exist.		
\item The set $U:=U_1\cup U_2$ with $U_1:=\big\{x\in\R^N:x_1\ge0\hbox{ and }x_2^2+\cdots+x_N^2\le1\big\}$ and $U_2:=\big\{x\in\R^N:x_1\ge0\hbox{ and }(x_2-\sqrt{x_1})^2+x_3^2+\cdots+x_N^2\le e^{-x_1^2}\big\}$, fulfills~\eqref{hyp:U} (hence~\eqref{ass0},~\eqref{ass},~\eqref{ass-cpt} and~\eqref{Hloc} hold), but it violates~\eqref{dUrho}, and~\eqref{dH} fails. Moreover, the first limit in~\eqref{Hausdorff} exists whereas the second one does not.
\item The set $U:=\big\{x\in\R^N:|x_N|\le e^{-|x'|^2}\big\}$ violates both~\eqref{hyp:U} and~\eqref{dUrho}, and~\eqref{ass0},~\eqref{ass},~\eqref{ass-cpt},~\eqref{Hloc} and~\eqref{dH} all fail, with $w(e)$ and $\W$ given by~\eqref{FGgeneral} and~\eqref{asspre}. Moreover, the two limits in~\eqref{Hausdorff} exist but do not coincide.
\item The set $U:=\big\{x\in\R^N:x_N\le \sqrt{|x'|}\big\}$ fulfills~\eqref{hyp:U} and~\eqref{dUrho} (hence \eqref{ass0}, \eqref{ass}, \eqref{ass-cpt}, \eqref{Hloc} and~\eqref{dH} all hold), but both limits in~\eqref{Hausdorff} do not exist and $d_{\mc{H}}(t^{-1}E_\lambda(t),\mc{W})=+\infty$ for all $\lambda\in(0,1)$ and $t>0$.	
\end{itemize}
The details about the above counter-examples can be found in \cite[Section 6]{HR1}.


\section{Flattening properties in the case of subgraphs}\label{secepi}
	
In this section, we focus on the important class of initial conditions which are characteristic functions of subgraphs in $\R^N$. Up to rotation, let us consider graphs in the direction~$x_N$, hence initial conditions $u_0$ given by
\be\label{defu0bis}
u_0(x',x_N)=\begin{cases} 1 & \text{if }x_N\le\gamma(x'),\\
0 & \text{otherwise},
\end{cases}
\ee
that is, $u_0=\1_U$ with $U$ given by~\eqref{Ugamma} and $\gamma\in L^\infty_{loc}(\R^{N-1})$. From parabolic estimates, $u(t,x',x_N)\to0$ as~$x_N\to+\infty$ and~$u(t,x',x_N)\to1$ as $x_N\to-\infty$, locally uniformly in $(t,x')\in[0,+\infty)\times\R^{N-1}$. Furthermore,~$\partial_{x_N}u<0$ in $(0,+\infty)\times\R^N$, by the strong parabolic maximum principle. As~a consequence, for every $t>0$, $x'\in\R^{N-1}$ and~$\lambda\in(0,1)$, there exists a unique $X_\lambda(t,x')\in\R$ such that
\be\label{defX}
u(t,x',X_\lambda(t,x'))=\lambda,
\ee
and the function $(\lambda,t,x')\mapsto X_\lambda(t,x')$ is actually continuous in $(0,1)\times(0,+\infty)\times\R^{N-1}$. In other words, the sets $E_\lambda(t)$ given in~\eqref{defElambda} are the open subgraphs of $x'\mapsto X_\lambda(t,x')$. Theorems~\ref{th1}-\ref{th4} applied to this case give some information on the shape of the graphs of $X_\lambda(t,\cdot)$ at large time and large space in terms of~$\gamma$, provided the assumptions of these theorems are fulfilled (see the previous section). We are now interested in the {\em local-in-space} behavior of the graphs of $X_\lambda(t,\cdot)$ at large time. Let us first point out that, because of the asymmetry of the roles of the steady states $0$ and $1$ (assuming Hypothesis~\ref{hyp:minimalspeed}), the behavior of the graphs of $X_\lambda(t,\cdot)$ will be radically different depending on the profile of the function~$\gamma$ at infinity. Consider the particular case~$\gamma(x')=\alpha\,|x'|$. Whatever $\alpha$ may be, the graphs of the functions $X_\lambda(t,\cdot)$ look like the sets $\{x\in\R^N:\dist(x,U)=c^*t\}$ at large time~$t$, by Theorem~\ref{th4}. For each~$t>0$, the set $\{x\in\R^N:\dist(x,U)=c^*t\}$ is a shift of $\partial U$ in the direction~$x_N$ and therefore has a vertex if $\alpha>0$, whereas it is $C^1$ if $\alpha\le0$. Of course, for each~$t>0$, in both cases~$\alpha>0$ and~$\alpha\le0$, each level set of $u$ (that is, each graph of $X_\lambda(t,\cdot)$) is of class~$C^2$ from the implicit function theorem and the negativity of~$\partial_{x_N}u$. Nevertheless, the previous observations imply that there should be a difference between the flattening properties of the level sets of $u$ according to the coercivity of~$\gamma$ at~infinity.

The following result deals with the non-coercive case, i.e., $\limsup_{|x'|\to+\infty}\gamma(x')/|x'|\le0$.

\begin{theorem}\label{thm:subgraph}{\rm{\cite{HR3}}}
Assume that Hypothesis~$\ref{hyp:minimalspeed}$ holds. Let~$u$ be the solution of~\eqref{homo} with an initial datum~$u_0$ given by~\eqref{defu0bis}. If
\Fi{gamma<0}
\limsup_{|x'|\to+\infty}\frac{\gamma(x')}{|x'|}\leq0,
\Ff
then, for every $\lambda\in[\theta,1)$ and every basis $(\mathrm{e}'_1,\cdots,\mathrm{e}'_{N-1})$ of $\R^{N-1}$, there holds
\be\label{liminf}
\liminf_{t\to+\infty}\Big(\min_{|x'|\le R,\,1\le i\le N-1}|\nabla_{\!x'}X_\lambda(t,x')\cdot\mathrm{e}'_i|\Big)\ \longrightarrow\ 0\ \hbox{ as }R\to+\infty,
\ee
where $\theta\in(0,1)$ is given by Proposition~$\ref{pro:spreadingFL}$ and Hypothesis~$\ref{hyp:invasion}$, and $X_\lambda$ by~\eqref{defX}. In particular, if~$N=2$, $\liminf_{t\to+\infty}\big(\min_{[-R,R]}|\partial_{x_1}X_\lambda(t,\cdot)|\big)\to0$ as $R\to+\infty$, for every $\lambda\in[\theta,1)$.
\end{theorem}

Roughly speaking, Theorem~\ref{thm:subgraph} says that the level set of any value $\lambda\in[\theta,1)$ becomes almost flat in some directions along sequences of points of $\R^N$ and along sequences of times diverging to $+\infty$. We point out that the estimates on $\nabla_{\!x'}X_\lambda(t,x')$ immediately imply analogous ones on $\nabla_{\!x'}u(t,x',X_\lambda(t,x'))$, because
$$\nabla_{\!x'}u(t,x',X_\lambda(t,x'))=-\partial_{x_N}u(t,x',X_\lambda(t,x'))\nabla_{\!x'}X_\lambda(t,x')$$
and $\partial_{x_N}u$ is bounded in $[1,+\infty)\times\R^N$ by parabolic estimates. Hence,~\eqref{liminf} implies that
$$\liminf_{t\to+\infty}\Big(\min_{|x'|\le R,\,1\le i\le N-1}|\nabla_{\!x'}u(t,x',X_\lambda(t,x'))\cdot\mathrm{e}'_i|\Big)\ \longrightarrow\ 0\ \hbox{ as }R\to+\infty.$$
for every $\lambda\in[\theta,1)$ and every basis $(\mathrm{e}'_1,\cdots,\mathrm{e}'_{N-1})$ of $\R^{N-1}$. The proof of~\eqref{liminf} is done in~\cite{HR3} by way of contradiction and uses the fact that $\lambda$ is larger than or equal to the quantity $\theta\in(0,1)$ given by Proposition~\ref{pro:spreadingFL} and Hypothesis~\ref{hyp:invasion}. If $f>0$ in~$(0,1)$, Hypothesis~\ref{hyp:invasion} is satisfied for any $\theta\in(0,1)$  and thus the
conclusion~\eqref{liminf} of Theorem~\ref{thm:subgraph} holds for any~$\lambda\in(0,1)$. Furthermore, if $f$ satisfies~\eqref{HTconditions}, then~\eqref{liminf} can be strengthened as $\liminf_{t\to+\infty}\big(\min_{|x'|\leq R}|\nabla_{\!x'}X_\lambda(t,x')|\big)\to0$ as $R\to+\infty$ for every $\lambda\in(0,1)$, see~\cite{HR3}.

Without~\eqref{gamma<0}, property~\eqref{liminf} does not hold in general  (immediate couterexamples are solutions whose level sets are parallel hyperplanes which are not orthogonal to~$\mathrm{e}_N$). Moreover, if one assumes $\liminf_{|x'|\to+\infty}\gamma(x')/|x'|\ge0$ instead of~\eqref{gamma<0}, the conclusion~\eqref{liminf} does not hold either in general (counter-examples are given by rotated bistable $V$-shaped fronts, from~\cite{HMR1,HMR2,NT,RR1}). However, with~\eqref{gamma<0}, we expect that the liminf of the min can be replaced by a limit in~\eqref{liminf}, without reference to the size $R$, namely, we propose the following.

\begin{conjecture}\label{conj1}
Under the assumptions of Theorem~$\ref{thm:subgraph}$, the conclusion~\eqref{liminf} can be strengthened 
by the limit, for every $\lambda\in[\theta,1)$,
\be\label{limit2}
\nabla_{\!x'}X_\lambda(t,x')\longrightarrow\ 0\ \hbox{ as }\text{$t\to+\infty$, \ locally uniformly in $x'\!\in\!\R^{N-1}$}.
\ee
\end{conjecture}

We emphasize that, even under the assumption~\eqref{gamma<0}, property~\eqref{limit2} does not hold in general {\em uniformly} with respect to $x'\in\R^{N-1}$ (for instance, in dimension $N=2$,~\eqref{limit2} fails if $\gamma'=1$ in $(-\infty,-1]$ and $\gamma'=-1$ in $[1,+\infty)$).

On the other hand, a strong support to the validity of Conjecture~\ref{conj1} is provided by the conclusion of Theorem~\ref{th4}. Indeed, it asserts that, for any $\lambda\in(0,1)$, $E_\lambda(t)\sim U+B_{c^*t}$ for $t$ large, in the sense of~\eqref{dH}, and one can check that condition~\eqref{gamma<0} entails that the exterior unit normals to the set $U+B_{c^*t}$ at the points $(x',x_N)\in\partial(U+B_{c^*t})$ (whenever they exist) approach the vertical direction $\mathrm{e}_N=(0,\cdots,0,1)$ as $t\to+\infty$, locally uniformly with respect to $x'\in\R^{N-1}$. Hence the same is expected to hold for the sets $E_\lambda(t)$, which is what~\eqref{limit2} asserts. This kind of arguments can be made rigorous, leading to a result which is a weaker statement than Conjecture~\ref{conj1}, that is, $\min_{|x'|\le\beta t}|\nabla_{x'}X_\lambda(t,x')|\to0$ as $t\to+\infty$ for every $\beta>0$, see~\cite{HR3}. A weaker statement than Conjecture~\ref{conj1} is also derived in~\cite{HR3} in the case where $f$ satisfies Fisher-KPP condition~\eqref{fkpp} below, namely, there holds that $\liminf_{t\to+\infty}\big(\max_{|x'|\le A}|\nabla_{x'}X_\lambda(t,x')|\big)=0$ for every $\lambda\in(0,1)$ and $A>0$ (see also Theorem~\ref{thm:DGe} below).

As for the full Conjecture~\ref{conj1}, the following result shows that, under assumption~\eqref{gamma<0}, the conclusion~\eqref{limit2} holds for initial data~$u_0$ having an asymptotically conical support, or being asymptotically $x'$-spherically-symmetric and nonincreasing. Notice that the following result uses the weaker Hypothesis~\ref{hyp:invasion} instead of Hypothesis~\ref{hyp:minimalspeed}.
	
\begin{theorem}\label{thm:conical}{\rm{\cite{HR3}}}
Assume that Hypothesis~$\ref{hyp:invasion}$ holds. Let $u$ be the solution of~\eqref{homo} with an initial datum $u_0$ given by~\eqref{defu0bis}, where $\gamma$ satisfies one of the following assumptions:
\begin{itemize}
\item[{\rm{(i)}}] either $\gamma$ is of class~$C^1$ outside a compact set and there is $\ell\ge0$ such that
\be\label{hypconical}
\hspace{-8pt}\left\{\baa{ll}
\!\!\!\gamma'(x_1)\to\mp\ell\ \hbox{ as }x_1\to\pm\infty & \hbox{if $\!N\!=\!2$},\vspace{3pt}\\
\displaystyle\!\!\!\nabla\gamma(x')=-\ell\,\frac{x'}{|x'|}+O(|x'|^{-1-\eta})\ \hbox{ as }|x'|\to+\infty,\hbox{ for some $\eta>0$,} & \hbox{if }\!N\!\ge\!3;\eaa\right.
\ee
\item[{\rm{(ii)}}] or $\gamma$ is continuous outside a compact set and $\gamma(x')/|x'|\to-\infty$ as~$|x'|\to+\infty$;
\item[{\rm{(iii)}}] or $\gamma(x')=\Gamma(|x'-x'_0|)$ outside a compact set,  for some $x'_0\in\R^{N-1}$ and some con\-tinuous nonincreasing function $\Gamma:\R^+\to\R$;
\item[{\rm{(iv)}}]  or $\gamma(x')=\Gamma(|x'-x'_0|)$ outside a compact set,  for some $x'_0\in\R^{N-1}$ and some $C^1$ function $\Gamma:\R^+\to\R$ such that  $\Gamma'(r)\to0$ as $r\to+\infty$.
\end{itemize}
Then, for every $\lambda_0\in(0,1)$, there holds that
\be\label{limit3chi}
\nabla_{\!x'}X_\lambda(t,x')\longrightarrow\ 0
\hbox{ \ as $t\to+\infty$, locally in $x'\!\in\!\R^{N-1}$ and uniformly in $\lambda\!\in\!(0,\lambda_0]$}
\ee
and moreover
$$\nabla_{\!x'}u(t,x',x_N)\longrightarrow\ 0
\hbox{ \ as $t\to+\infty$, locally in $x'\!\in\!\R^{N-1}$ and uniformly in $x_N\!\in\!\R$}.$$
\end{theorem}

In dimension $N=3$, by writing $\gamma(x')=\tilde{\gamma}(r,\vartheta)$ in the standard polar coordinates, condition~\eqref{hypconical} means that $\partial_r\tilde{\gamma}(r,\vartheta)=-\ell+O(r^{-1-\eta})$ and $\partial_\vartheta\tilde{\gamma}(r,\vartheta)=O(r^{-\eta})$ as $r\to+\infty$.

It is easy to see that, even under Hypothesis~\ref{hyp:minimalspeed} (which is stronger than Hypothesis~\ref{hyp:invasion}), if~\eqref{hypconical} holds with $\ell>0$, then the convergence in~\eqref{limit3chi} cannot be uniform with respect to~$x'\in\R^{N-1}$. In other words, if the initial interface between the states $0$ and $1$ has a non-zero slope at infinity, then the level sets cannot become uniformly flat at large time. This observation naturally leads to the following conjecture.
	
\begin{conjecture}\label{conj:order1}
Assume that Hypothesis $\ref{hyp:minimalspeed}$ holds. Let $u$ be the solution of~\eqref{homo} with an initial datum $u_0$ given by~\eqref{defu0bis}. If
\Fi{gamma'to0}
\lim_{|x'|\to+\infty}\nabla\gamma(x')=0,
\Ff
then, for every $\lambda_0\in(0,1)$,
\be\label{Dchi0}
\nabla_{\!x'}X_\lambda(t,x')\longrightarrow\ 0\hbox{ \ as $t\to+\infty$, uniformly in $x'\!\in\!\R^{N-1}$ and in $\lambda\!\in\!(0,\lambda_0]$}
\ee
and moreover
\be\label{Du0}
\nabla_{\!x'}u(t,x)\longrightarrow\ 0\hbox{ \ as $t\to+\infty$, uniformly in $x\!\in\!\R^N$}.
\ee
\end{conjecture}
	
Properties~\eqref{Dchi0}-\eqref{Du0} obviously hold if $\gamma$ is constant. Furthermore, if condition~\eqref{gamma'to0} is replaced by the boundedness of $\gamma$, then, at least for some classes of functions $f$, pro\-perties~\eqref{Dchi0} (with $\lambda\in(0,\lambda_0]$ replaced by $\lambda\in[a,b]$, for some fixed $0<a\le b<1$) and~\eqref{Du0} hold: more precisely, if the function $f$ is of the bistable type~\eqref{bistable}, these properties follow from some results in~\cite{BH1,FM}, and the same conclusions hold for more general functions $f$ of the multistable type~\cite{P2} or for KPP type functions~$f$ satisfying~\eqref{fkpp} below or slightly weaker conditions, see~\cite{BH1,B,HNRR,L,U1}. Further estimates on the exact position of the level sets~$X_\lambda$ in the bistable or KPP cases have been established in~\cite{MN,MNT,RR2}. However, by considering some functions $\gamma$ with large local oscillations at infinity, it turns out that both conclusions of Conjecture~\ref{conj:order1} cannot hold if~\eqref{gamma'to0} is replaced by the weaker condition $\lim_{|x'|\to+\infty}\gamma(x')/|x'|=0$, see~\cite{HR3}.
	
To complete this section, let us point out that, under the assumptions of Theorems~\ref{th1} and~\ref{thm:subgraph}, the solution $u$ of~\eqref{homo} with~\eqref{defu0bis} propagates with speed~$c^*$ in the direction $\mathrm{e}_N=(0,\cdots,0,1)$, owing to~Theorem~\ref{th1}, that is, $w(\mathrm{e}_N)=c^*$ in~\eqref{ass}-\eqref{FGgeneral}. We conjecture that the solution $u$ then locally converges along its level sets to the front profile $\varphi$ with minimal speed $c^*$.

\begin{conjecture}\label{conj3}
Under the assumptions of Theorems~$\ref{th1}$ and $\ref{thm:subgraph}$, it holds, for every $\lambda\in(0,1)$, for every sequence $(t_n)_{n\in\N}$ diverging to $+\infty$, and for every bounded sequence $(x'_n)_{n\in\N}$ in~$\R^{N-1}$,
\be\label{convphi}
u(t_n+t,x'_n+x',X_\lambda(t_n,x'_n)+x_N)\longrightarrow\vp(x_N-c^*t+\vp^{-1}(\lambda))\quad\text{as }n\to+\infty,
\ee
in $C^{1;2}_{loc}(\R_t\times\R_{x'}^{N-1})$ and uniformly with respect to $x_N\in\R$. If one further assumes~\eqref{gamma'to0}, then the above limit holds for every sequence $(x'_n)_{n\in\N}$ in $\R^{N-1}$, bounded or not.
\end{conjecture}

By~\cite{HR3}, the second conclusion does not hold in general if assumption~\eqref{gamma'to0} is replaced by $\lim_{|x'|\to+\infty}\gamma(x')/|x'|=0$. On the other hand, Conjecture~\ref{conj3}, and especially its second part, holds if $\gamma$ is bounded, for some classes of functions $f$, see~\cite{BH1,MN,MNT,P2,RR2}.
	

\section{Asymptotic one-dimensional symmetry}\label{secaos}

Let us now present some results about the asymptotic one-dimensional symmetry, 
related to the Question~\ref{q2} presented in Section~\ref{sec:2questions}. 
They concern Fisher-KPP~\cite{F,KPP} functions $f$, that is, satisfying
\be\label{fkpp}
f(0)\!=\!f(1)\!=\!0,\ f(s)\!>\!0\hbox{ for all }s\!\in\!(0,1),\hbox{ and }s\mapsto\frac{f(s)}{s}\hbox{ is nonincreasing in $(0,1]$}.
\ee
In this case the {\em hair trigger effect} holds, i.e., Hypothesis~\ref{hyp:invasion} is fulfilled for any~$\theta,\rho\!>\!0$, moreover Hypothesis~\ref{hyp:minimalspeed} also holds and the minimal speed with the properties stated in Proposition~\ref{pro:spreadingFL} is explicit: $c^*=2\sqrt{f'(0)}$, see~\cite{AW,KPP}.

\begin{theorem}\label{thm:DG}{\rm{\cite{HR2}}}
Assume that $f$ is of the Fisher-KPP type~\eqref{fkpp}. Let $u$ be the solution of~\eqref{homo} with an initial datum $u_0=\1_{U}$ such that $U\subset\R^N$ satisfies 
\Fi{UUdelta}
\exists\,\delta>0,\quad d_{\mc{H}}(U,U_\delta)<+\infty.
\Ff 
Assume moreover that $U$ is convex, or more generally, that there is a convex set~$U'$ such that $d_{\mc{H}}(U,U')<+\infty$. Then, any function $\psi\in\Omega(u)$, with $\Omega(u)$ defined by~\eqref{def:Omega}, is of the form $\psi(x)=\Psi(x\.e)$, for some~$e\in\Sph$ and a function $\Psi:\R\to\R$ which is either constant or strictly monotone.
\end{theorem}

Theorem~\ref{thm:DG} extends the asymptotic one-dimensional symmetry property known to hold when $U$ is bounded, as a consequence of some results of Jones~\cite{J} valid for even more general functions $f$ provided $U$ contains a large ball, as well as when $U$ is the subgraph of a bounded function, by~\cite{BH1,B,HNRR,L,U1}. Conversely, the asymptotic one-dimensional symmetry fails when $U$ is ``$V$-shaped'', i.e., the union of two non-parallel half-spaces, which fulfills~\eqref{UUdelta} but is not at a finite Hausdorff distance from a convex set nor it satisfies~\eqref{ballcone}. For such an initial datum, the $\O$-limit set of the solution contains elements which are not one-dimensional.
	
Condition~\eqref{UUdelta} means that there exists $R>0$ such that, for any~$x\in U$, there is a ball $B_\delta(x_0)\subset U$ with $|x-x_0|<R$. It is fulfilled in particular if~$U$ satisfies a uniform interior ball condition. One can show that, in dimension $N=2$, for a convex set $U$, property~\eqref{UUdelta} is equivalent to require that $U$ has nonempty interior. The role of~\eqref{UUdelta} is cutting off regions of $U$ playing a negligible role in the large-time behavior of the solutions. This assumption is necessary, otherwise one could consider a $V$-shaped set $\tilde U$ and then take $U:=\tilde U\,\cup\,\bigcup_{k\in\Z^N}B_{e^{-|k|^2}}(k)$, which is at finite Hausdorff distance from the convex set~$\R^N$ but does not satisfy~\eqref{UUdelta}, and the associated solution violates the asymptotic one-dimensional symmetry (because it essentially behaves at large time as the solution associated with $\tilde{U}$) see~\cite{HR2}. The aforementioned two examples show that the answer to Question~\ref{q2} cannot be positive without any assumption on~$U$.

The idea of the proof of \thm{DG} in~\cite{HR2} consists in reducing to a case where it is possible to apply the reflection argument of Jones~\cite{J}, which is valid for more general functions $f$ but fails when $U$ is unbounded. This is achieved by an approximation of the solution through a suitable truncation of its initial support. In order to control the error, new types of supersolutions initially supported in exterior domains are used, which are also employed in the proofs of the results of Section~\ref{sec:FG}.

As a matter of fact, the convex-proximity assumption on $U$ in \thm{DG} is a very special case of a geometric hypothesis under which the one-dimensional symmetry holds. Namely, for a given nonempty set $U\subset\R^N$ and a given point $x\in\R^N$, we let $\pi_x$ denote the set of orthogonal projections of~$x$ onto $\ol U$, i.e.,
$$\pi_x:=\big\{\xi\in\ol U: |x-\xi|=\dist(x,U)\big\},$$
and, for $x\notin\ol U$, we define the {\em opening} function as follows:
$$\mc{O}(x):=\sup_{\xi\in\pi_x,\,y\in U\setminus\{\xi\}}\frac{x-\xi}{|x-\xi|}\.\frac{y-\xi}{|y-\xi|},$$
with the convention that $\mc{O}(x)=-\infty$ if $U=\emptyset$ or $U$ is a singleton (otherwise $-1\leq\mc{O}(x)\leq1$).
When $\mc{O}(x)\neq-\infty$, one has $\mc{O}(x)=\cos\alpha$, where $\alpha$ is the infimum among all~$\xi\in\pi_x$ of half the opening of the largest exterior cone to $U$ at~$\xi$ having axis~$x-\xi$. Here is our most general asymptotic symmetry result.

\begin{theorem}\label{thm:DGgeneral}{\rm{\cite{HR2}}}
Assume that $f$ is of the Fisher-KPP type~\eqref{fkpp}. Let $u$ be a solution of~\eqref{homo} with an initial datum $u_0=\1_{U}$ such that $U\subset\R^N$ satisfies~\eqref{UUdelta} and moreover
\Fi{ballcone}
\lim_{R\to+\infty}\bigg(\,\sup_{x\in\R^N,\,\dist(x,U)=R}\mc{O}(x)\bigg)\leq 0.
\Ff		
Then any function in $\Omega(u)$ is one-dimensional and, in addition, it is either constant or strictly monotone, in the sense of Theorem~$\ref{thm:DG}$.
\end{theorem}

It is understood that the left-hand side in condition~\eqref{ballcone} is equal to $-\infty$ (hence the condition is fulfilled) if $\sup_{x\in\R^N}\dist(x,U)<+\infty$ 
(and indeed in such a case one has by~\eqref{UUdelta} that $\sup_{x\in\R^N}\dist(x,U_\delta)<+\infty$ for some 
$\delta>0$, which implies that  $u(t,x)\to1$ uniformly in $x\in\R^N$ as~$t\to+\infty$ due to
 Proposition~\ref{pro:spreadingFL}, hence in particular the asymptotic one-dimensional symmetry holds). 
 We also point out that the limit in~\eqref{ballcone} always exists, because the involved quantity is nonincreasing with respect to~$R$, see~\cite{HR2}. Hypothesis~\eqref{ballcone} means that the angle $\alpha$ in the definition of $\mc{O}(x)$ tends to a value larger than or equal to~$\pi/2$ (i.e., the exterior cone contains a half-space) as $\dist(x,U)\to+\infty$. Theorem~\ref{thm:DGgeneral} yields Theorem~\ref{thm:DG} because, firstly, convex sets satisfy $\mc{O}(x)\leq0$ for every $x\notin\ol U$ (actually, they are characterized by such condition in the class of closed sets) and, secondly, if \eqref{ballcone} holds for a given set, then it holds true for any set at finite Hausdorff distance from it, by~\cite{HR2}. However, the class of sets satisfying~\eqref{ballcone} is wider. It contains for instance the subgraphs of functions with {\em vanishing global mean},~i.e., $U$ of the type~\eqref{Ugamma} with $\gamma\in L^\infty_{loc}(\R^{N-1})$ such~that
\Fi{VGM}
\frac{\gamma(x')-\gamma(y')}{|x'-y'|}\longrightarrow 0\as|x'-y'|\to+\infty
\Ff
(see Corollary~\ref{cor:DGVGM} below for the precise result in this case). 
Actually, with~\eqref{Ugamma} and~\eqref{VGM}, any $\psi\in\Omega(u)$ is of the form $\psi(x)=\Psi(x_N)$, 

Theorems~\ref{thm:DG} and~\ref{thm:DGgeneral} are concerned with {\it locally uniform} convergence properties along sequences of times $(t_n)_{n\in\N}$ diverging to $+\infty$ and sequences of points $(x_n)_{n\in\N}$. As a matter of fact, a {\it uniform} asymptotic property can be derived, from the conclusions of Theorems~\ref{thm:DG} or~\ref{thm:DGgeneral}. It is expressed in terms of the eigenvalues of the Hessian matrices $D^2u(t,x)$ (with respect to the $x$ variables). For a symmetric real-valued matrix $A$ of size $N\times N$, let $\lambda_1(A)\le\cdots\le\lambda_N(A)$ denote its eigenvalues, and let
$$\sigma_k(A):=\sum_{1\le j_1<\cdots<j_k\le N}\lambda_{j_1}(A)\times\cdots\times\lambda_{j_k}(A),\ \ \ \ 1\le k\le N,$$
be the elementary symmetric polynomials of eigenvalues of $A$ ($\sigma_k(D^2u(t,x))$ is also called $k$-Hessian).

\begin{theorem}\label{thm:global}{\rm{\cite{HR2}}}
Let $f$ and $u$ be as in Theorem~$\ref{thm:DGgeneral}$. Then, for every $2\le k\le N$, $\sigma_k(D^2u(t,x))\to0$ as $t\to+\infty$ uniformly in $x\in\R^N$.
\end{theorem}

Theorem~\ref{thm:global} is proved using the asymptotic local one-dimensional symmetry given by Theorem~\ref{thm:DGgeneral}, together with standard parabolic estimates. We point out that, if $\psi:\R^N\to\R$ is of class $C^2$ and one-dimensional, then $\sigma_k(D^2\psi(x))=0$ for all $2\le k\le N$ and $x\in\R^N$, since the quantities $\sigma_k(D^2\psi(x))$ involve sums of products of at least two eigenvalues of $D^2\psi(x)$ (but $\sigma_1(D^2\psi(x))\neq0$ in general). However, the converse property is immediately not true (for instance, the function $\psi:(x_1,x_2)\mapsto x_1^2+x_2$ satisfies $\sigma_2(D^2\psi(x_1,x_2))=0$ for all $(x_1,x_2)\in\R^2$, but it is not one-dimensional).

Once the asymptotic one-dimensional symmetry and monotonicity properties are established, it is natural to ask what are the directions in which the solution actually becomes locally one-dimensional. Namely, we investigate the~set
\Fi{E}
\begin{split}
\mc{E}:=\big\{\, & e\in\Sph\, :\, \exists\,\psi\in\O(u)\text{ such that $\,\psi(x)\equiv\Psi(x\.e)$}\\
& \text{for some strictly decreasing function $\Psi\in C^2(\R)$}\big\}.
\end{split}
\Ff
Under the assumptions of Theorems~\ref{thm:DG} or~\ref{thm:DGgeneral}, the set $\mc{E}$ is then the set of the directions of  decreasing monotonicity of all non-constant elements of $\Omega(u)$ (by the direction of decreasing monotonicity of a --necessarily one-dimensional by Theorems~\ref{thm:DG} or~\ref{thm:DGgeneral}-- non-constant function $\psi\in\Omega(u)$, we mean the unique $e\in\Sph$ such that $\psi(x)=\Psi(x\cdot e)$ for all $x\in\R^N$, with $\Psi$ decreasing). The constant functions $\psi$ are excluded in the above definition, this is necessary because they are one-dimensional in every direction. Thus, a direction~$e$ belongs to~$\mc{E}$ only if, along some diverging sequences of times, the solution flattens in the directions orthogonal to~$e$ but not in the direction $e$, along some sequence of points. We characterize the set $\mc{E}$ in terms of the initial support $U$.

\begin{theorem}\label{thm:E}{\rm{\cite{HR2}}}
Let $f$ and $u$ be as in Theorem~$\ref{thm:DGgeneral}$. Then $\mc{E}$ defined in~\eqref{E} is given by
\[
\begin{split}
\mc{E}=\Big\{ \, & e\in\Sph\ :\  \displaystyle\frac{x_n-\xi_n}{|x_n-\xi_n|}\to e\ \hbox{as $n\to+\infty$, for some sequences $(x_n)_{n\in\N}$, $\seq{\xi}$ in $\R^N$}\\
&\text{such that $\,\dist(x_n,U)\to+\infty$ as $n\to+\infty$ and $\xi_n\in \pi_{x_n}$ for all $n\in\N$}\Big\}.
\end{split}
\]
In particular, $\mc{E}=\emptyset$ if and only if $U$ is relatively dense in~$\R^N$ or $U=\emptyset$. Moreover, for any $e\in\mc{E}$, any sequence $(x_n)_{n\in\N}$ in $\R^N$, and any sequence $(t_n)_{n\in\N}$ diverging to $+\infty$ such that $u(t_n,x_n+x)\to\Psi(x\cdot e)$ as $n\to+\infty$ locally uniformly in~$x\in\R^N$, with $\Psi$ strictly decreasing, one necessarily has $\dist(x_n,U)\sim c^*t_n$ as $n\to+\infty$.
\end{theorem}

We remark that, without the assumption~\eqref{UUdelta}, the last statement of \thm{E} may immediately fail. Indeed, if $U=\{0\}$ then $u(t,x)\equiv0$ for all $t>0$, $x\in\R^N$, hence $\mc{E}=\emptyset$, but~$U\neq\emptyset$ is not relatively dense in $\R^N$.

When $U$ is bounded with non-empty interior, it follows from Theorem~\ref{thm:E} that $\mc{E}=\Sph$. On the one hand, this conclusion gives an additional property --namely the strict monotonicity-- with respect to the asymptotic symmetry result contained in~\cite{J}. On the other hand, still when $U$ is bounded, the same conclusion is also a consequence of~\cite{D,RRR}, where it is proved by a completely different argument than in~\cite{HR2}. The characterization of the directions of asymptotic strict monotonicity in the case of unbounded sets $U$ is more involved. The proof of Theorem~\ref{thm:E} is based on an argument by contradiction and on the acceleration of the solutions when they become less and less steep.

\thm{E} implies that if~$U$~is of class~$C^1$ then~$\mc{E}$ is contained in the closure of the set of the outward unit normal vectors~to~$U$. If $U$ is convex then $\mc{E}$ coincides with the closure of the set of outward unit normal vectors to all half-spaces containing $U$. When $U$ is the subgraph of a function $\gamma$ with vanishing global mean, i.e.~satisfying~\eqref{VGM}, then it turns out that~$\mc{E}=\{\mathrm{e}_N\}$, with $\mathrm{e}_N=(0,\cdots,0,1)$, as the following result shows.

\begin{corollary}\label{cor:DGVGM}{\rm{\cite{HR2}}}
Assume that $f$ is of the Fisher-KPP type~\eqref{fkpp}. Let $u$ be the solution of~\eqref{homo} with an initial datum $u_0=\1_{U}$, where $U$ is given by~\eqref{Ugamma} with $\gamma\in L^\infty_{loc}(\R^{N-1})$ satisfying~\eqref{VGM}. Then any function $\psi\in\Omega(u)$ is of the form $\psi(x',x_N)\equiv\Psi(x_N)$ for all $(x',x_N)\in\R^{N-1}\times\R$, with $\Psi$ either constant or strictly decreasing. Moreover, $\mc{E}=\{\mathrm{e}_N\}$.
\end{corollary}

Since, by parabolic estimates, the convergence in the definition~\eqref{def:Omega} of the $\O$-limit set holds true in $C^2_{loc}(\R^N)$ up to subsequences, Corollary~\ref{cor:DGVGM} implies that
$$\nabla_{\!x'} u(t,x',x_N)\to0\ \hbox{ as }t\to+\infty,\ \text{ uniformly with respect to $(x',x_N)\in\R^{N-1}\times\R$}.$$
A way to interpret this result is that the oscillations of the initial datum are ``damped'' as time goes~on through some kind of averaging process. We point out that Corollary~\ref{cor:DGVGM} does not imply the existence of a function $\Psi:\R^+\times\R\to\R$ such that $u(t,x',x_N)-\Psi(t,x_N)\to0$ as $t\to+\infty$ uniformly in $(x',x_N)\in\R^{N-1}\times\R$, and indeed such a function $\Psi$ does not exist in general (as shown in~\cite{RR2} when $N=2$ and the limits $\lim_{x'\to\pm\infty}\gamma(x')$ exist but do not coincide). Condition~\eqref{VGM} is satisfied in particular when $\gamma$ is bounded, and in such a case the conclusion of Corollary~\ref{cor:DGVGM} can also be deduced from~\cite{BH1,B,HNRR,L,U1}.

It is possible to relax the uniform mean condition~\eqref{VGM} of $\gamma$ in Corollary~\ref{cor:DGVGM}, at the price of restricting the $\O$-limit set. Similarly, the arguments of the proof of \thm{DGgeneral} can somehow be localized. Loosely speaking, if one focuses on the asymptotic one-dimensional property around a given direction, the global geometric assumption~\eqref{ballcone} can be restricted to the points $x$ around that direction, and hypothesis~\eqref{UUdelta} can be relaxed too. This leads us to introduce the notion of {\it $\O$-limit set in a direction $e\in\Sph$} of a solution $u$, defined as
$$\begin{array}{ll}
\Omega_e(u):=\big\{&\!\!\!\!\psi\in L^\infty(\R^N)\ :\   u(t_n,x_n+\.)\to\psi\text{ in $L^\infty_{loc}(\R^N)$}\\
& \!\!\!\! \text{for some sequences $(t_n)_{n\in\N}$ in $\R^+$ diverging to $+\infty$},\\
& \!\!\!\! \text{and $(x_n)_{n\in\N}$ in $\R^N\setminus\{0\}$ such that $x_n/|x_n|\to e$ as $n\to+\infty$}\,\big\}\ \subset\ \Omega(u).\end{array}$$

\begin{theorem}\label{thm:DGe}{\rm{\cite{HR2}}}
Assume that $f$ is of the Fisher-KPP type~\eqref{fkpp}. Let $u$ be the solution of~\eqref{homo} with an initial condition  $u_0\!=\!\1_U$, where $U$ has nonempty interior and is such that
\Fi{Ucontained}
U\subset\big\{(x',x_N)\in\R^{N-1}\times\R \ :\ x_N\leq\gamma(x')\big\},
\Ff
for a function $\gamma\in L^\infty_{loc}(\R^{N-1})$ satisfying~\eqref{gamma<0}. Then, any function $\psi\in\Omega_{\mathrm{e}_N}(u)$ satisfies $\psi(x',x_N)\equiv\Psi(x_N)$ in~$\R^N$, 
with $\Psi$ either constant or strictly decreasing. In particular, $\nabla_{\!x'}u(t,x',x_N)\to0$ as $t\to+\infty$ locally in $x'\!\in\!\R^{N-1}$ and uniformly in $x_N\!\in\![R,+\infty)$, for any $R\in\R$. Moreover if the inclusion is replaced by an equality in~\eqref{Ucontained}, then $\nabla_{\!x'}u(t,x',x_N)\to0$ as $t\to+\infty$ locally in $x'\!\in\!\R^{N-1}$ and uniformly in $x_N\!\in\!\R$.
\end{theorem}

To complete this section, we propose a list of open questions and conjectures related to the previous results. First of all, under the Fisher-KPP assumption~\eqref{fkpp}, let $\varphi$ be the traveling front profile with minimal speed, that is, for each $e\in\Sph$, $\varphi(x\cdot e-c^*t)$ satisfies~\eqref{homo} with $0=\varphi(+\infty)<\varphi<\varphi(-\infty)=1$ and $c^*=2\sqrt{f'(0)}$. Based on Theorems~\ref{thm:DGgeneral} and~\ref{thm:E}, and together with the definition~\eqref{E} of $\mc{E}$, we propose the following.

\begin{conjecture}
Let $f$ and $u$ be as in Theorem~$\ref{thm:DGgeneral}$. Then,
\[
\O(u)=\big\{\, 0,\ 1,\ x\mapsto\vp(x\.e +a)\ : \ e\in\mc{E},\ a\in\R\,\big\}.
\]
\end{conjecture}

The above conjecture is known to be true in the case where $U$ is bounded with non-empty interior, by~\cite{D,RRR}, 
as well as when $U$ is the subgraph of a bounded function, or more generally when there are two half-spaces $H$ and $H'$ --necessarily with parallel boundaries-- such that $H\subset U\subset H'$, by~\cite{BH1,B,HNRR,L,U1}.

We remark that the assumption~\eqref{ballcone} of Theorem~\ref{thm:DGgeneral} is stable by bounded perturbations of the sets $U$. We could then wonder whether the \aos\ is also stable with respect to bounded perturbations of the initial support. Namely, if the solution to~\eqref{homo} with an initial datum $\1_{U}$ satisfying~\eqref{UUdelta} is asymptotically locally planar, and if $U'\subset\R^N$ satisfies~\eqref{UUdelta} and $d_{\mc{H}}(U',U)<+\infty$, then is the solution to~\eqref{homo} with initial datum $\1_{U'}$ asymptotically locally planar as well?

One can also wonder whether the reciprocal of Theorem~\ref{thm:DGgeneral} is true, in the following sense: if the \aos\ holds for a solution $u$ of~\eqref{homo} with initial datum $\1_{U}$ and $U$ satisfying~\eqref{UUdelta}, does necessarily $U$ fulfill~\eqref{ballcone}? The answer is immediately seen to be negative in general: take for instance $U$ given by $U=\bigcup_{n\in\N}[2^n,2^n+1]\times\R^{N-1}$, which fulfills~\eqref{UUdelta} but not~\eqref{ballcone}, while $u$ --~hence any element of $\Omega(u)$~-- is one-dimensional, depending on the variable $x_1$ only. However, the question is open if~$U$ is connected.

The results of this section are concerned with the Fisher-KPP equation, i.e.~when~$f$ satisfies~\eqref{fkpp}. However, the same question about the \aos\ can be asked for more general reaction terms $f$, still with $f(0)=f(1)=0$ and satisfying Hypothesis~\ref{hyp:minimalspeed}, or simply the invasion property stated in Hypothesis~\ref{hyp:invasion}. First of all, the condition~\eqref{UUdelta} should be strengthened, by requiring $\delta$ to be larger than the quantity $\rho$ in Hypothesis~\ref{hyp:invasion} (indeed, if $f$ is for instance of the bistable type~\eqref{bistable} with $\int_0^1f(s)ds>0$, then by~\cite{DM1,Z} there is $\delta_0>0$ such that the solution to~\eqref{homo} with initial condition $u_0=\1_{B_{\delta_0}}$ converges uniformly as $t\to+\infty$ to a ground state, that is, a positive radial solution converging to $0$ as $|x|\to+\infty$, hence $u$ is not asymptotically locally planar). For general functions $f$ for which Hypothesis~\ref{hyp:invasion} holds, if $U$ is bounded and $U_{\rho}\neq\emptyset$, then the solutions to~\eqref{homo} with initial condition $\1_U$ are known to be asymptotically locally planar, by~\cite{J}. The same conclusion holds for bistable functions $f$ of the type~\eqref{bistable} if there are two half-spaces $H$ and $H'$ --necessarily with parallel boundaries-- such that $H\subset U\subset H'$, by~\cite{BH1,FM,MN,MNT} (see also~\cite{P2} for the case of more general functions $f$). On the other hand, still for bistable functions $f$ of the type~\eqref{bistable} for which $\int_0^1f(s)ds>0$, the solutions $u$ to~\eqref{homo} with initial condition $\1_U$ are not asymptotically locally planar if $U$ is V-shaped, that is, if it is the union of two half-spaces with non-parallel boundaries, by~\cite{HMR1,HMR2,NT,RR1}. These known results lead us to formulate the following De Giorgi type conjecture for the solutions of the reaction-diffusion equation~\eqref{homo} beyond the Fisher-KPP case. 

\begin{conjecture}\label{conj:DG}
Assume that the invasion property, namely Hypothesis~$\ref{hyp:invasion}$, holds for some $\rho>0$. Let $u$ be the solution to~\eqref{homo} with an initial datum $u_0=\1_{U}$ such that $U\subset\R^N$ satisfies~$d_{\mc{H}}(U,U_\rho)<+\infty$ and~\eqref{ballcone}. Then any function in $\Omega(u)$ is one-dimensional and, in addition, it is either constant or strictly monotone.
\end{conjecture}

We point out that, when $u_0$ is given by~\eqref{defu0bis} with $\gamma$ satisfying~\eqref{gamma'to0}, Conjecture~\ref{conj:order1} would imply the validity of the asymptotic one-dimensional symmetry.

Let us also mention another natural question related to the preservation of the convexity of the upper level sets of $u$ when $u_0=\1_U$ and $U$ is convex. It is known from~\cite{BL,IST} that, if $U$ is convex, then the solution of the heat equation $\partial_tu=\Delta u$ is quasi-concave at each $t>0$, that is, for each $t>0$ and $\lambda\in\R$, the upper level set $\{x\in\R^N:u(t,x)>\lambda\}$ is convex. The same conclusion holds for~\eqref{homo} set in bounded convex domains instead of $\R^N$, and under some additional assumptions on $f$, by~\cite{IS}. A natural question is to wonder for which class of functions $f$ this property still holds for~\eqref{homo} in $\R^N$.

Notice finally that, for any solution $u$ to~\eqref{homo}, for any sequence $(t_n)_{n\in\N}$ diverging to~$+\infty$, and for any sequence $(x_n)_{n\in\N}$ in $\R^N$, the functions $u(t_n+\cdot,x_n+\cdot)$ converge locally uniformly in $\R\times\R^N$, up to extraction of a subsequence, to an entire solution to~\eqref{homo} (that~is, a solution for all~$t\in\R$). Remembering Theorem~\ref{thm:DG} on the \aos\ for the solutions to~\eqref{homo} with $u_0=\1_U$ and $U$ convex, and having in mind the question of the previous paragraph on the convexity of the upper level sets, it is then natural to ask the following: if an entire solution $v:\R\times\R^N\to[0,1]$ to~\eqref{homo} is quasi-concave for every $t\in\R$, is $v(t,\cdot)$ necessarily one-dimensional for every $t\in\R$?


\section{The logarithmic lag in the KPP case}\label{sec:KPP}	

Assume in this section that $f$ satisfies the Fisher-KPP condition~\eqref{fkpp}. We recall that Hypotheses~\ref{hyp:invasion} and~\ref{hyp:minimalspeed} are fulfilled, and the minimal speed $c^*$ of traveling fronts connec\-ting~$1$ to~$0$ is given by $c^*=2\sqrt{f'(0)}$. It is known from~\cite{B,HNRR,L,NRR,U1} that, in the one-dimensional case, the solution $u$ of~\eqref{homo} with initial condition $u_0=\1_{\R^-}$ is such that
$$\sup_{x\in\R}\Big|u(t,x)-\varphi\Big(x-c^*t+\frac{3}{c^*}\,\ln t+x_0\Big)\Big|\to0\hbox{ as $t\to+\infty$},$$
for some $x_0\in\R$. Hence, there is a lag by $(3/c^*)\ln t$ of the position of the level sets of~$u$ behind the position $c^*t$ given by the spreading speed. In dimension $N=2$, for initial conditions trapped between two shifts of $\1_{\R\times\R^-}$, there is a bounded function $a$ such that
$$\sup_{(x_1,x_2)\in\R^2}\Big|u(t,x_1,x_2)-\varphi\Big(x_2-c^*t+\frac{3}{c^*}\,\ln t+a(t,x_1)\Big)\Big|\to0\hbox{ as $t\to+\infty$},$$
see~\cite{RR2}. In any dimension $N\ge2$, if the nonnegative initial condition $0\not\equiv u_0\le1$ is compactly supported, there is a Lipschitz continuous function $a:\Sph\to\R$ such that
$$\sup_{x\in\R^N\setminus\{0\}}\Big|u(t,x)-\varphi\Big(|x|-c^*t+\frac{N+2}{c^*}\,\ln t+a\Big(\frac{x}{|x|}\Big)\Big)\Big|\to0\ \hbox{ as }t\to+\infty,$$
see~\cite{D,G,RRR}. Notice that~$N+2=3+(N-1)$ corresponds to an additional lag by $((N-1)/c^*)\ln t$, compared with the $1$-dimensional case, which is due to the curvature of the level sets inherited from the fact that the initial condition is compactly supported.

Let us now consider the case of a solution to~\eqref{homo} with an initial condition given by~\eqref{defu0bis}
and investigate the lag between the position of the level sets of~$u$ behind $c^*t$ in the direction~$x_N$. 
Assuming that $\gamma$ is bounded from above, one infers by comparison that, up to an additive constant, the lag is between $(3/c^*)\ln t$ (the lag in the $1$-dimensional case) and $((N+2)/c^*)\ln t$ (the lag in the case of compactly supported initial conditions): namely, for every~$\lambda\in(0,1)$ and $x'\in\R^{N-1}$, under the notations~\eqref{defX}, the lag $c^*t-X_\lambda(t,x')$ satisfies
\be\label{lagbetween}
\frac{3}{c^*}\,\ln t+O(1)\leq c^*t-X_\lambda(t,x')\leq \frac{N+2}{c^*}\,\ln t+O(1)\ \ \hbox{ as }t\to+\infty.
\ee
But it is not clear in principle whether or not this lag is equal to one of these bounds or whether it takes intermediate values. The main result of this section states that the actual lag coincides with the upper bound
in~\eqref{lagbetween} provided that $\gamma$ tends to $-\infty$ at infinity faster than a logarithm with a 
suitable negative coefficient. Thus, in such a case, 
the position of the level sets of $u$ in the direction $x_N$ is the same as when the initial condition is compactly supported. Here is the precise result.

\begin{theorem}\label{thm:normdelay}{\rm{\cite{HR3}}}
Assume that $f$ is of the Fisher-KPP type~\eqref{fkpp} and let $u$ be the solution of~\eqref{homo} with an initial condition $u_0$ satisfying~\eqref{defu0bis}. If
\Fi{gamma<}
\limsup_{|x'|\to+\infty}\,\frac{\gamma(x')}{\ln(|x'|)}<-\frac{2(N-1)}{c^*},
\Ff
then
\be\label{lag1}
X_\lambda(t,x')=c^*t-\frac{N+2}{c^*}\,\ln t+O(1)\ \hbox{ as }t\to+\infty,
\ee
locally uniformly with respect to $\lambda\in(0,1)$ and $x'\in\R^{N-1}$, and the inequality ``$\leq$'' holds true in the above formula locally uniformly in $\lambda\in(0,1)$ and uniformly in $x'\in\R^{N-1}$.
\end{theorem}

If the upper bound for $\gamma$ in~\eqref{gamma<} is relaxed, we expect the lag of the solution with respect to the critical front to differ from the one associated with compactly supported initial data, that is $((N+2)/c^*)\ln t$. We derive the following lower bound for the lag.

\begin{proposition}\label{pro:lag}{\rm{\cite{HR3}}}
Assume that $f$ is of the Fisher-KPP type~\eqref{fkpp} and let $u$ be the solution of~\eqref{homo} with an initial condition $u_0$ satisfying~\eqref{defu0bis}. If there is $\sigma\ge-(N-1)$ such~that
\be\label{hyplag}
\limsup_{|x'|\to+\infty}\frac{\gamma(x')}{\ln|x'|}\le\frac{2\sigma}{c^*},
\ee
then, for any $\lambda\in(0,1)$,
\be\label{ineqlag}
X_\lambda(t,x')\leq c^*t-\frac{3-\sigma}{c^*}\,\ln t+o(\ln t)\as t\to+\infty,
\ee
locally uniformly with respect to $x'\in\R^{N-1}$.
\end{proposition}

Property~\eqref{ineqlag} means that the lag $c^*t-X_\lambda(t,x')$ is at least~$((3-\sigma)/c^*)\ln t+o(\ln t)$ as $t\to+\infty$. Notice that this holds even for positive $\sigma$. We conjecture that, if the limsup is replaced by a limit in~\eqref{hyplag} and the inequality by an equality, then the estimate~\eqref{ineqlag} should be sharp, namely, the lag should be
$$c^*t-X_\lambda(t,x')=\frac{3-\sigma}{c^*}\ln t+o(\ln t)\quad\text{as $t\to+\infty$},$$
for every $\lambda\in(0,1)$ and~$x'\in\R^{N-1}$. We emphasize that when $\sigma=0$, this formula would be coherent with the $1$-dimensional lag. This formula would also mean that the constant $-2(N-1)/c^*$ in~\eqref{gamma<} would be optimal for the lag to be equivalent to that of solutions with compactly supported initial conditions. Lastly, it would provide a continuum of lags with logarithmic factors ranging in the whole half-line $(-\infty,(N+2)/c^*]$. In particular, solutions with initial conditions of the type~\eqref{defu0bis} with~$\gamma(x')\sim(6/c^*)\ln|x'|$ as~$|x'|\to+\infty$ would have no logarithmic lag, i.e., the same position~$c^*t$ along the $x_N$-axis as the planar front moving in the direction $\mathrm{e}_N$, up to a $o(\ln t)$ term as $t\to+\infty$. Furthermore, if~$\gamma(x')\sim \kappa\ln|x'|$ as~$|x'|\to+\infty$ for some $\kappa>(6/c^*)$, then the logarithmic lag would be negative, i.e., the position of the solution would be ahead of that of the front by a logarithmic-in-time term (observe that the term is linear in time when $\gamma(x')\sim \alpha |x'|$ as~$|x'|\to+\infty$ with $\alpha>0$,  by Theorem~\ref{th1} and~\eqref{ass}-\eqref{FGgeneral}).

Theorem~\ref{thm:normdelay} allows us to prove part of Conjecture~\ref{conj1} about the flattening of the level sets under the hypotheses of that theorem.

\begin{corollary}\label{cor:DGlag}{\rm{\cite{HR3}}}
Assume that $f$ is of the Fisher-KPP type~\eqref{fkpp} and let $u$ be the solution of~\eqref{homo} with an initial condition $u_0$ satisfying~\eqref{defu0bis} and~\eqref{gamma<}. Then the following hold:
\begin{itemize}
\item[{\rm{(i)}}] the conclusion~\eqref{limit2} of Conjecture~$\ref{conj1}$ holds, and even locally in $\lambda\in(0,1)$, that is, $\nabla_{\!x'}X_\lambda(t,x')\to0$ as $t\to+\infty$, locally uniformly in $x'\in\R^{N-1}$ and $\lambda\in(0,1)$;
\item[{\rm{(ii)}}] for any $\lambda\in(0,1)$ and~$x'_0\in\R^{N-1}$, the function
$$\t u(t,x',x_N):=\lim_{s\to+\infty}u(s+t,x',X_\lambda(s,x'_0)+x_N),$$
which exists $($up to subsequences$)$ locally uniformly in $(t,x',x_N)\in\R\times\R^N$, is independent of $x'$ and satisfies $\lim_{x_N\to-\infty}\t u(t,x_N+c^*t)=1$ and $\lim_{x_N\to+\infty}\t u(t,x_N+c^*t)=0$, uniformly with respect to $t\in\R$.
\end{itemize}
\end{corollary}

Corollary~\ref{cor:DGlag} shows that, as $t\to+\infty$, the solution approaches a one-dimensional entire solution whose level sets move in the direction $\mathrm{e}_N$ with average velo\-city equal to~$c^*$. It is then natural to expect that $\t u(t,x_N)=\varphi(x_N-c^*t+\varphi^{-1}(\lambda))$ for all $(t,x_N)\in\R^2$, where $\varphi$ is the front connecting~$1$ and~$0$ with minimal speed $c^*$. That would correspond to property~\eqref{convphi} in Conjecture~\ref{conj3}. By comparison and some arguments based on the number of intersections of solutions to~\eqref{homo} in dimension $1$, it can be shown that $\t u(t,x_N)\ge\varphi(x_N-c^*t+\zeta)$ in $\R^2$, for some $\zeta\in\R$. But the proof of~\eqref{convphi} would still require additional arguments.


\end{document}